\documentclass{amsart}
\usepackage[table]{xcolor} 
\usepackage[francais,english]{babel}
\usepackage[colorlinks,citecolor=darkblue,urlcolor=blue,bookmarks=false,hypertexnames=true]{hyperref} 
\usepackage{amsrefs}
\usepackage{amsmath}
\usepackage{amsfonts}
\usepackage{graphicx,color}
\usepackage{amssymb}
\newtheorem{theo}{Theorem}[section]

\newtheorem{prop}{Proposition}[section]

\newtheorem{lem}{Lemma}[section]

\newtheorem{step}{Step}[section]
\definecolor{darkblue}{HTML}{004C93} 
\definecolor{MainRed}{rgb}{.6, .1, .1}
\def\E {E}
\def \H {\tilde{H}_0(\mathbb{R}_+^N)}
\def \rr {\mathbb{R}}
\def \crit{2_s^{\star}}
\def \crito{2_0^{\star}}

\def \vv {\mathbb{\Vert}}

\def \e {\mathbb{R}_+^{N}}
\def \rn {\mathbb{R}^N}

\def \nn {\mathbb{N}}

\def \la {\mathbb{\lambda}}

\def \ep {\epsil\nu_{\alpha}on}

\def \bb {\hbox}

\def \gg {\gamma }

\def \ep {\epsilon}

\date{}

\title[Fourth order Hardy-Sobolev equations]{Fourth order Hardy-Sobolev equations: Singularity and doubly critical exponent}

\author{Hussein Cheikh Ali}
\address{D\'epartement de Math\'ematiques, Universit\'e libre de Bruxelles, CP 214, Boulevard du Triomphe, B-1050 Bruxelles, Belgium}
\email{Hussein.cheikh-ali@ulb.ac.be}
\begin{document}
	\maketitle
	\pdfbookmark[0]{\contentsname}{Contents}
	\begin{abstract}
In dimension $N\geq 5$, and for $0< s<4$ with $\gamma\in\rr$, we study the existence of nontrivial weak solutions for the doubly critical problem 
		$$\Delta^2 u-\frac{\gamma}{|x|^4}u= |u|^{2^\star_0-2}u+\frac{|u|^{\crit-2}u}{|x|^s}\hbox{ in }\rr_+^N,\; u=\Delta u=0\hbox{ on }\partial \rr_+^N,$$
where  $\crit:=\frac{2(N-s)}{N-4}$ is the critical Hardy–Sobolev exponent. For $N\geq 8$ and $0< \gg<\frac{(N^2-4)^2}{16}$,  we show   the existence of nontrivial solution using the Mountain-Pass theorem by Ambrosetti-Rabinowitz. The method used is based on the existence of extremals for certain Hardy-Sobolev embeddings that we prove in this paper.
	\end{abstract}
\section{Introduction}
Let  $\rr_+^N:=\{x\in \rr^N\, | \,x_1>0\}$ be the half-space domain, $N\geq 5$ and $0\leq s<4$. In this work, we establish the
existence of nontrivial weak solutions for the following doubly critical problem:
\begin{equation}\label{eq:doublehardybi}
	\left\{\begin{array}{ll}
		\Delta^2 u -\frac{\gg}{|x|^4} u=|u|^{2^\star_0-2}u+\frac{|u|^{\crit-2}u}{|x|^s} &\hbox{ in } \rr_+^N, \\
		u=\Delta u = 0  &\hbox{ on } \partial \rr_+^N,
	\end{array}\right.
\end{equation}
where $\Delta=\hbox{div}\left( \nabla\right)$ is the Laplacian, $\gg\in \rr$ and $\crit:=\frac{2(N-s)}{N-4}$ is the critical Hardy–Sobolev exponent. In order to write a variational formulation, the relevant space is the following: given $\Omega$ be a smooth domain in $\rr^N$, define 
\begin{equation*}
\tilde{H}_0(\Omega)=\hbox{ completion of }\{u\in C^2_c(\overline{\Omega})\hbox{ s.t. }u_{|\partial\Omega}\equiv 0\}\hbox{ for the norm }u\mapsto \Vert\Delta u\Vert_2.
\end{equation*}
\noindent Note that $\tilde{H}_0(\rn)$ is the usual Beppo-Levi space $D^{2,2}(\rn)$. We say that $u\in \H$ is a weak solution of \eqref{eq:doublehardybi}, if  \begin{eqnarray*}
\int_{\e}  \Delta u \Delta \varphi  \, dx-\gg\int_{\e}\frac{u\varphi}{|x|^4}\, dx=\int_{\e} \left|u \right|^{\crito-2}u\varphi\, dx+\int_{\e}\frac{\left|u\right|^{\crit-2}u\varphi}{ |x|^s}\,dx,
\end{eqnarray*}
 for all $\varphi\in \H$. For any $u\in \H$, Sobolev's embedding yields $u\in L^{\crito}(\e)\cap L^{\crit}(\e, |x|^{-s})$ (see \eqref{ineq:sobo} below), and therefore the definition of weak solution makes sense.\par
 \medskip \noindent The existence of weak solution to  \eqref{eq:doublehardybi} on $\rr^N$ has been studied by Filippucci-Pucci-Robert \cite{FPR} for the p--Laplacian. Equations like \eqref{eq:doublehardybi} have been studied for the Fractional Laplacian (see Ghoussoub-Shakerian \cite{GS}) and  the bi-laplacian in $\rn$ (see  Bhakta \cites{B1,B2} and Bhakta-Musina \cite{BM}).\par 
 \smallskip\noindent In the present paper, we tackle this type of nonlinear singular problems on $\rr_+^N$ when $0\in \partial \rr_+^N$. When $0\in\partial\Omega$ where  $\Omega$ is a smooth domain of $\rr^N$, $(N\geq 3 )$, the existence of solutions for the corresponding 2nd order with only nonlinearity have been studied by  Ghoussoub–Kang \cite{GK} and studied by Chern-Lin \cite{CL} and Ghoussoub-Robert \cite{GRGAFA,GR}. For non-smooth domains modeled on cones, we refer to the more recent works of Cheikh-Ali \cite{HCA1,HCA2}.\par

\medskip\noindent Our main result is the following:
\begin{theo}\label{theo2}
	Let $N\geq 8$, $0<s<4$ and $0<\gg <\frac{(N^2-4)^2}{16}$. Then, there exists a nontrivial weak solution of \eqref{eq:doublehardybi}.
\end{theo}
Let us discuss the hypothesis of the theorem. Our problem depends of the value of the constant $\gg$.  From here, given an arbitrary domain $\Omega\subset \rn$,  $N\geq 5$,  we define the Hardy-Rellich constant
  \begin{equation}\label{def:hardy:X}
  	\gg_H(\Omega,X):=\inf\left\lbrace 	\frac{	\int_{\Omega}\left| \Delta u\right| ^2\, dx}{\int_{\Omega}\frac{u^2}{|x|^4}\, dx}; u\in X\backslash\{0\} \right\rbrace,
  \end{equation}
for a suitable space $X$.  There are several references around this constant. A first version of this Hardy inequality has been introduced by  Rellich in 1953 \cites{R1, R2} (see also Mitidieri \cite{EM}), and reads 
\begin{equation*}\label{ineq:hardy:rn}
\gg_{H}(\rr^N,C^2_c(\rn))= \frac{N^2(N-4)^2}{16}.
\end{equation*}
When $\Omega$ is a bounded domain in $\rr^N$ with $0\in \Omega$, Perez-Llamos-Primo \cite{PP} proved that
$$\gg_{H}(\Omega,H^{2}(\Omega)\cap H_0^1(\Omega))= \frac{N^2(N-4)^2}{16}.$$
In particular, the value of the constant is independent of the domain as long as $0$ is an interior point.

  \medskip \noindent The situation is different when $0\in\partial\Omega$. Let us consider cones. For any regular domain $\Sigma$  in the unit sphere $\mathbb{S}^{N-1}$, we define the cone 
 $$C_{\Sigma}:=\{r\sigma \, | \,r>0, \sigma \in  \Sigma\}.$$
Caldiroli-Musina \cite{CM} proved that 
$$\gg_{H}(C_{\Sigma},X_0(C_\Sigma)) =\hbox{dist}\left(-\frac{N(N-4)}{4}, \Lambda(\Sigma)\right)^2,$$
where,
$$X_0(\Omega):=\{u\in C^2(\overline{\Omega})\cap C^2_c(\rn\setminus\{0\})\hbox{ s.t. }u_{|\partial\Omega}\equiv 0\},$$
and $\Lambda(\Sigma)$ is the spectrum of the Laplace–Beltrami operator on $\Sigma$. We denote that, for instance:
  \begin{itemize}
  	\item If $\Sigma=\mathbb{S}^{N-1}$, then $C_{\Sigma}= \rr^N\backslash \{0\}$. Therefore, we have
  	\begin{align*}
  		\gg_{H}(\rr^N\backslash \{0\},X_0(\rr^N\backslash \{0\})) &=\min_{k\in\{0,1,...\}}\left|\frac{N(N-4)}{4}+k(N-2+k) \right|^2=\frac{N^2(N-4)^2}{16}.
  	\end{align*} 
  	\item If $\Sigma$ is the half-sphere $S^{N-1}_+$, then $C_{\Sigma}=\e$. Therefore, we have 
  	\begin{equation}\label{eq:hardyinC2}
  		\gg_{H}(\e,X_0(\e)) =\min_{k\in\{1,2,...\}}\left|\frac{N(N-4)}{4}+k(N-2+k) \right|^2=\frac{(N^2-4)^2}{16}.
  	\end{equation} 
  \end{itemize}
In the sequel, we write for convenience $$\gg_{H,+}:=\gg_{H}(\e,\H).$$
We will prove in  Section \ref{sec:HardyRellich} that
  $$\gg_{H,+} =\gg_{H}(\e,X_0(\e))=\frac{(N^2-4)^2}{16}.$$  
In order to tackle the nonlinear problem \eqref{eq:doublehardybi}, let us consider the Sobolev inequality 
\begin{equation}\label{ineq:sobo}
\left(\int_{\rn}\frac{|u|^{\crit}}{|x|^s}\, dx\right)^{\frac{2}{\crit}}\leq C\int_{\rn}|\Delta u|^2\, dx\hbox{ for all }u\in C^\infty_c(\rn).
\end{equation}
\noindent Let us fix a domain $\Omega\subset\rn$. Interpolating the Hardy inequality \eqref{def:hardy:X} and the Sobolev inequality \eqref{ineq:sobo} and using that $X_0(\Omega)\subset C^\infty_c(\rn)$, given $s\in [0,4]$, we get the existence of $C(N,s,\Omega)>0$ such that
 \begin{equation}\label{eq:hardys}
   	\left( \int_{\Omega}|x|^{-s}|u|^{\crit}\,dx\right)^{\frac{2}{\crit}} \leq C(N,s,\Omega) \int_{\Omega}\left| \Delta u\right|^2\, dx \bb{ for any } u\in C^\infty_c(\Omega).
	\end{equation}
Using again \eqref{def:hardy:X}, for any $s\in [0,4]$ and any $\gamma<\gamma_H(\Omega, \tilde{H}_0(\Omega))$, we get the Hardy-Sobolev inequality
 \begin{equation*}\label{ineq:hardysobolev:omega}
 \left( \int_{\Omega}|x|^{-s}|u|^{\crit}\,dx\right)^{\frac{2}{\crit}} \leq 	C\,  \int_{\Omega}\left( \left| \Delta u\right|^2-\gamma\frac{u^2}{|x|^4}\right) \, dx \bb{ for all }u\in \tilde{H}_0(\Omega).
 \end{equation*}
Let us define 
\begin{equation}\label{eq:qgammas}
	Q_{\gamma, s}(\Omega)=\inf\Big\{I_{\gamma, s}^{\Omega}(u)/u\in \tilde{H}_0(\Omega)\backslash\{0\}\Big\}, \hbox{ with }I_{\gamma, s}^{\Omega}(u):=\frac{\int_{\Omega}\left( \left| \Delta u\right| ^2 -\gamma\frac{u^2}{|x|^4}\right) \,dx}{\left(\int_{\Omega}\frac{|u|^{\crit}}{|x|^s}\, dx\right)^{\frac{2}{\crit}}}.
\end{equation}
In order to prove the existence of weak solutions for \eqref{eq:doublehardybi}, we need extremals for the best constant $Q_{\gg,s}(\e)$ in \eqref{eq:qgammas}, that is $u\in\H\setminus\{0\}$ such that $I_{\gamma, s}^{\e}(u)=Q_{\gg,s}(\e)$.  If $u\in \H \backslash\{0\}$ is an extremal for  $Q_{\gamma, s}(\e)>0$,  then, up to a constant factor, $u$ is a solution to the following Euler-Lagrange equation:
\begin{equation*}\label{eq:hardybi}
	\left\{\begin{array}{cc}
		\Delta^2 u -\frac{\gg}{|x|^4} u=\frac{|u|^{\crit-2}u}{|x|^s} &\hbox{ in } \rr_+^N, \\
		u=\Delta u = 0  &\hbox{ on } \partial \rr_+^N.
	\end{array}\right.
\end{equation*} 
The question of the existence of extremals is answered in the following theorem that we prove in Section \ref{sec:extr}:
\begin{theo}\label{theo1}
	For $s\in [0,4)$ and $\gg < \gg_{H,+}$, we have that 
	\begin{itemize}
		\item[(a)] If $\{s>0\}$ or $\{ s=0, \gg>0, N\geq 8\}$, then there are extremals for $Q_{\gg,s}(\rr_+^N)$.
		\item[(b)] If $\{s=0\bb{ and } \gg\leq 0\}$, there are no extremals for $Q_{\gg,s}(\rr_+^N)$.
		\item[(c)] If there are no extremals for $Q_{\gg,0}(\e)$, then $Q_{\gg,0}(\e)=S_N$,
	\end{itemize}
where 
\begin{equation}\label{eq:bestconstantsoblevrn}
S_N:=\inf_{u\in \tilde{H}_0(\rn)\backslash\{0\}}	\frac{\int_{\rn}\left| \Delta u\right| ^2\, dx}{\left(\int_{\rn}|u|^{\crito}\, dx \right)^{\frac{2}{\crito}}}.
\end{equation}
\end{theo}
\noindent The remaining case, that is $\{s=0\, ,\,\gg>0\hbox{, and }N=5,6,7\}$, is not clear. This limitation is due to the lack of localization of the $L^2-$norm in the computation \eqref{eq:i1ep}.

 \smallskip The study of fourth-order Hardy-Sobolev problem with a  singularity on the boundary of a smooth domain is the object of the work \cite{HCA3}.\par
 
 \medskip\noindent With all these elements, we will construct  weak solutions for the  doubly critical problem \eqref{eq:doublehardybi} by finding critical points of corresponding functional on $\H$. The method to obtain this critical points is via the Mountain-Pass theorem of Ambrosetti and Rabinowitz. Since \eqref{eq:doublehardybi} is invariant under the conformal one parameter transformation group:
$$\left\{\begin{array}{cccc}
	T_r: & \H  &\to & \H \\
	& u& \mapsto & T_r[u]:=\{x\mapsto r^{\frac{N-4}{2}}u(rx)\} 
\end{array}\right\}, \bb{ where } r>0,$$
then the Mountain-Pass theorem will not yield critical points, but only 
Palais-Smale sequences. We will use the strategy of Filippucci-Pucci-Robert \cite{FPR}. As in \cite{FPR}, the main difficulty will be the  asymptotic competition between the energy carried by the two critical nonlinearities. Hence, the crucial point here is to balance the competition to avoid the domination of one term over the other. Otherwise, there is vanishing of the weakest one, and we get solution for the same equation but with only one critical nonlinearity. To deal with this problem, we will choose a Palais-Smale sequence at a minimax energy level. In such a way, after a careful analysis of concentration, we will show that there is a balance between the energies of the two nonlinearities mentioned above, and therefore none can dominate the other. This will yield a solution to \eqref{eq:doublehardybi}.

\section{Profile of solutions and study the value of Hardy Rellich Constant $\gg_{H,+}$ }\label{sec:HardyRellich}
We first determine the value of the Hardy Rellich constant $\gg_{H,+}$ when $0\in \partial \rr_+^N$.  
	\begin{lem}\label{lemma:extraonhalfspace}
	We have that 
	\begin{eqnarray*}\label{eq:lemmaext}
	\gg_{H,+}:=\gg_{H}\left( \e,\H\right)=\gg_{H}\left( \e,X_0(\rn_+)\right)=\frac{(N^2-4)^2}{16}.
	\end{eqnarray*}
\end{lem}
\medskip\noindent \textit{Proof of Lemma \ref{lemma:extraonhalfspace}:} We prove the lemma in two steps proving each an inequality. 

\begin{step}\label{step1:lemma1}
	We claim that 
	\begin{align}\label{eq:inversHS}
	\int_{\rr_+^N}\left| \Delta v\right|^2\, dx \geq 	\frac{(N^2-4)^2}{16}\, \int_{\rr_+^N}\frac{v^2}{|x|^4}\, dx  \bb{ for all } v\in C^{2}_c(\overline{\rr_+^N})\hbox{ s.t. }v_{|\partial\rn_+}=0.
	\end{align}
\end{step}
\noindent{\it Proof of Step \ref{step1:lemma1}:} Fix $v\in C^{2}_c(\overline{\rr_+^N})$ such that $v_{\partial\rn_+}=0$. We choose $\eta\in C^\infty(\rr^N)$ 
 such that $\eta(x)_{|B_1(0)}\equiv 0$, $\eta(x)_{|B_2(0)^c}\equiv 1$  and $0\leq \eta\leq 1$.
For $\ep>0$, we set $v_\ep(x):=\eta_\ep(x) v(x) , \bb{ where } \eta_\ep(x):=\eta(\frac{x}{\ep}) \bb{ for all } x\in\rr^N.$ We have that $v_\ep\in X_0(\rn_+)$. It follows from Caldiroli-Musina \cite{CM} (see \eqref{eq:hardyinC2}) that 
\begin{align}\label{eq:inversHSep}
\int_{\rr_+^N}\left|\Delta v_\ep\right|^2\, dx	 \geq \frac{(N^2-4)^2}{16}\, \int_{\rr_+^N}\frac{v_\ep^2}{|x|^4}\, dx.
\end{align}
\medskip\noindent\textit{Step \ref{step1:lemma1}.1} We claim that 
\begin{align}\label{eq:delta1ep}
	\int_{\rr_+^N}\left| \Delta v_\ep \right|^2\,dx=	\int_{\rr_+^N} \left|\eta_\ep \right|^2 \left|\Delta v\right|^2\,dx+o(1) \bb{ as }\ep\to 0.
\end{align}
\medskip\noindent\textit{Proof of this claim:}  For convenience we define $A_{+,\ep}:=\rr_+^N\cap \left( B_{2\ep}(0)\backslash B_{\ep}(0) \right).$
Using the definition of $v_\ep$ yields
\begin{equation}\label{eq:variationDelta}
	\int_{\rr_+^N}|\Delta v_\ep|^2\,dx=\int_{\rr_+^N}|\eta_\ep|^2|\Delta v|^2 \,dx+R_{\ep},
\end{equation}
with
\begin{eqnarray}\label{aDelta0ep}
	R_\ep&:=& \int_{A_{+,\ep} }\Big[  |\Delta \eta_\ep|^2v^2+4\nabla \eta_\ep\cdot\nabla v \Delta(\eta_\ep)v\\
	&&+\,2\, \Delta(\eta_\ep)\eta_\ep v\Delta v+4\, (\nabla \eta_\ep\cdot\nabla v)^2+4\, \nabla \eta_\ep\cdot\nabla v \eta_\ep\Delta v\Big]\,dx. \nonumber
\end{eqnarray}
We claim that
\begin{equation}\label{eq:Rep}
	R_\ep=O\left(\int_{A_{+,\ep}}\frac{|v|^2}{|x|^4}\,dx+\left( \int_{A_{+,\ep}}\frac{|v|^2}{|x|^4}\,dx\right)^{\frac{1}{2}}+\left( \int_{A_{+,\ep}}\frac{|v|^2}{|x|^4}\,dx\right)^{\frac{1}{4}}	\right).
\end{equation}
\begin{proof} We estimate each term of $R_\ep$. Since $v_{|\partial\rn_+}=0$, integrating by parts yield
	\begin{eqnarray*}
		\int_{A_{+,\ep}} \nabla \eta_\ep\cdot\nabla v \Delta(\eta_\ep)v\, dx
		&=&O	\left( -\int_{A_{+,\ep}}\left(  |\Delta \eta_\ep|^2+ \nabla \eta_\ep\cdot \nabla (\Delta \eta_\ep)\right) v^2 \,dx\right. \\
		&&\left. +\int_{\rr_+^N\cap\partial\left(B_{2\ep}(0)\backslash B_{\ep}(0) \right) }v^2 \Delta \eta_\ep\partial_{\nu}\eta_\ep \,d\sigma\right) ,
	\end{eqnarray*}
	where $\nu$ is the outer normal vector of $B_{2\ep}(0)\backslash B_{\ep}(0)$. Since $\partial_{\nu} \eta_\ep=0$ on $\rr_+^N\cap\partial\left(B_{2\ep}(0)\backslash B_{\ep}(0) \right)$, we have
	\begin{align}
		\int_{A_{+,\ep}} \nabla \eta_\ep\cdot\nabla v  \Delta(\eta_\ep)v\, dx
		&=O\left(  \frac{1}{\ep^4}\vv (\Delta \eta)^2+ \nabla \eta\cdot \nabla (\Delta \eta)\vv_{\infty}\int_{A_{+,\ep}} v^2\, dx\right)\nonumber \\
		&=O\left( \int_{A_{+,\ep}}\frac{|v|^2}{|x|^4}\, dx\right).\label{e12ep}  
	\end{align}
By Hölder's inequality and $v\in C^2_c(\overline{\e})$, we get 
\begin{eqnarray}
	\int_{A_{+,\ep}} \Delta(\eta_\ep)\eta_\ep v\Delta v\, dx&&=O\left( \left( \int_{A_{+,\ep}}\left|\Delta v \right|^2\, dx\right)^{\frac{1}{2}}\left( \int_{A_{+,\ep}}\left|\Delta(\eta_\ep)\eta_\ep v \right|^2\, dx\right)^{\frac{1}{2}}\right)  \nonumber\\
	&&=O\left(  \left( \int_{A_{+,\ep}}\frac{|v|^2}{|x|^4}\,dx\right)^{\frac{1}{2}}	\right). \label{e14ep} 
\end{eqnarray}
It follows from the Cauchy-Schwarz and H\"older inequalities and integrations by parts that
	\begin{eqnarray}
		\int_{A_{+,\ep}}\Big( \nabla \eta_\ep\cdot\nabla v\Big)^2 \,dx &=&O\left( \int_{A_{+,\ep}}|\nabla \eta_\ep|^2|\nabla v|^2 \, dx\right)\nonumber\\
		&=&O\left( \int_{A_{+,\ep}} v\Delta v  |\nabla \eta_\ep|^2 \, dx+\int_{A_{+,\ep}}v\nabla v\nabla(|\nabla \eta_\ep|^2 )\, dx\right)\nonumber\\
		&=&O\left( \left( \int_{A_{+,\ep}}\frac{|v|^2}{|x|^4}\,dx\right)^{\frac{1}{2}}+\int_{A_{+,\ep}}\frac{|v|^2}{|x|^4}\, dx\right).\label{e13ep}
	\end{eqnarray}
Using again Hölder's inequality yields
\begin{eqnarray}
	\int_{A_{+,\ep}} \nabla \eta_\ep\cdot \nabla v \eta_\ep\Delta v\,dx &=&O\left( \left( \int_{A_{+,\ep}}|\nabla \eta_\ep|^2|\nabla v|^2 \, dx\right)^{\frac{1}{2}}\right)\nonumber\\
	&=&O\left( \left( \int_{A_{+,\ep}}\frac{|v|^2}{|x|^4}\,dx\right)^{\frac{1}{2}}+\left( \int_{A_{+,\ep}}\frac{|v|^2}{|x|^4}\,dx\right)^{\frac{1}{4}}\right).\label{e15ep}
\end{eqnarray}
We inject \eqref{e12ep}, \eqref{e14ep}, \eqref{e13ep} and \eqref{e15ep}  in \eqref{aDelta0ep}, we obtain \eqref{eq:Rep}. 
\end{proof}
\noindent It follows from \eqref{eq:Rep} that $R_\ep=o(1)$ as $\ep \to 0$. Therefore, by \eqref{eq:variationDelta}
we obtain as $\ep \to 0$ that  \eqref{eq:delta1ep}. This ends the proof of Step \ref{step1:lemma1}.\qed\par  
\medskip \noindent Using again the inequality \eqref{eq:inversHSep} and by \eqref{eq:delta1ep}, we find that
\begin{align*}
	\int_{\rr_+^N}|\eta_\ep|^2\left| \Delta v\right|^2\, dx+o(1)\geq \frac{(N^2-4)^2}{16}\, \int_{\rr_+^N}\frac{|\eta_\ep|^2 v^2}{|x|^4}\, dx .
\end{align*}
Therefore, passing $\ep \to 0$, we get \eqref{eq:inversHS}.  This proves Step \ref{step1:lemma1}.\qed 
\begin{step}\label{step2:lemma1} We claim that
\begin{align}\label{eq:inversHS1}
	 \int_{\rr_+^N}|\Delta u|^2\, dx\geq\frac{(N^2-4)^2}{16}\, \int_{\rr_+^N}\frac{u^2}{|x|^4}\, dx  \bb{ for all } u\in \H.
\end{align} 
\end{step}
\noindent{\it Proof of Step \ref{step2:lemma1}:} We fix $u\in  \H$. We then take a sequence $(u_n)_{n}$ such that $u_n\in C^{2}_c(\overline{\rr_+^N})$ and $u_n(x)=0$ for all $x\in\partial\rn_+$ and $\lim\limits_{n\to +\infty}u_n=u$ for the norm $\Vert\Delta\cdot\Vert_2$. Therefore, 
\begin{equation*}
	\lim_{n\to +\infty} \int_{\rr_+^N}\frac{u_n^2}{|x|^4}\, dx=\int_{\rr_+^N}\frac{u^2}{|x|^4}\, dx\hbox{ and }	\lim_{n\to +\infty} \int_{\rr_+^N}\left|\Delta u_n\right|^2\, dx =\int_{\rr_+^N} \left| \Delta u\right|^2\, dx.
\end{equation*}
It then follows from \eqref{eq:inversHS} that 
\begin{align*}
	\int_{\rr_+^N}|\Delta u_n|^2\, dx\geq \frac{(N^2-4)^2}{16}\, \int_{\rr_+^N}\frac{u_n^2}{|x|^4}\, dx .
\end{align*}
Letting  $n\to+\infty$ and we get \eqref{eq:inversHS1}. This proves Step \ref{step2:lemma1}. \qed \par 
\begin{step}\label{step3:lemma1} We claim that \eqref{eq:inversHS1} is optimal.
\end{step}
\noindent{\it Proof of Step \ref{step3:lemma1}:} This will be achieved via test-function estimates. We  define 
$$v(x):=x_1|x|^{-\frac{N-2}{2}}\hbox{ for all } x\in\rr_+^N\backslash \{0\}.$$
We have that $v\in C^2(\overline{\e})$, $v=0$ on $\partial\e$ and 
$$-\Delta v(x)=\frac{N^2-4}{16}\frac{v(x)}{|x|^{2}}.$$
Let $\varphi,\psi\in C^\infty(\rr^n)$ such that  
$$\left\{\begin{array}{l}
\varphi(0)=0,\, |\varphi(x)|\leq c|x|\hbox{ if }|x|<1\hbox{ and }\varphi(x)=1\hbox{ if }|x|\geq 1;\\
\psi(x)=1\hbox{ if }|x|< 1 \hbox{ and } \psi(x)=0 \hbox{ if }|x|>2,
\end{array}\right.$$
for some constant $c>0$. For $0<\ep\ll 1$, we define the function $v_{\ep}\in \H$ as follows:
\begin{equation*}
v_{\ep}(x):=\varphi\left(\frac{x}{\ep}\right)\psi(\ep x)v(x)=\left\{\begin{array}{cc}
\varphi\left(\frac{x}{\ep}\right) v(x) &\hbox{ if }|x|<\ep,\\
 v(x) &\hbox{ if }\ep\leq |x|<\frac{1}{\ep},\\
 \psi(\ep x)v(x) &\hbox{ if }|x|\geq\frac{1}{\ep}.
\end{array}\right.
\end{equation*}
It follows from the definition of $v_\ep$ that
\begin{equation}\label{hb1}
	\int_{\rn_+\backslash \overline{B}_{\ep^{-1}}(0)}\left| \Delta v_{\ep}\right| ^2dx=O(1) \bb{ and } 	\int_{\rn_+\cap B_{\ep}(0)}\left| \Delta v_{\ep}\right| ^2\, dx=O(1) \bb{ as } \ep\to 0.
\end{equation}
It remains just one factor to calculate, we have 
\begin{eqnarray}
	\int_{B_{\ep^{-1}}(0)\backslash\overline{B}_{\ep}(0)}\left| \Delta v_{\ep}\right| ^2\, dx
	&=&\frac{(N^2-4)^2}{16}\int_{B_{\ep^{-1}}(0)\backslash\overline{B}_{\ep}(0)}\frac{v_{\ep}^2}{|x|^4}\, dx\label{hb30}.
\end{eqnarray} 
Moreover, we calculate that
\begin{eqnarray}
	\int_{B_{\ep^{-1}}(0)\backslash\overline{B}_{\ep}(0)}\frac{v_{\ep}^2}{|x|^4}dx&=&\int_{B_{\ep^{-1}}(0)\backslash\overline{B}_{\ep}(0)}x_1^2|x|^{-N-2}dx\nonumber\\
	&=&2w(2)\ln\left(\frac{1}{\ep}\right)\label{hb40},
\end{eqnarray}
where $w(2):=\int_{\mathbb{S}_+^{N-1}}x_1^2\, d\sigma$.
It follows from \eqref{hb1}, \eqref{hb30} and \eqref{hb40} that
\begin{eqnarray}
	\int_{\rr_+^N}\left| \Delta v_{\ep}\right| ^2\, dx=2w(2)\frac{(N^2-4)^2}{16}\ln\left(\frac{1}{\ep}\right)+O(1),\label{hb50}
\end{eqnarray}
when $\ep \to 0$. Using again the definition of $v_\ep$ yields,
\begin{eqnarray}\label{hb05}
	\int_{\rn_+\backslash \overline{B}_{\ep^{-1}}(0)}\frac{v_{\ep}^2}{|x|^4}\, dx=O(1) &\bb{ and }& \int_{ \rn_+\cap B_{\ep}(0)}\frac{v_{\ep}^2}{|x|^4}\, dx=O(1).
\end{eqnarray}
Therefore, it follows from \eqref{hb40}, we get as $\ep\to 0$ that
\begin{eqnarray}
	\int_{\rr_+^N}\frac{v_{\ep}^2}{|x|^4}\, dx=2w(2)\ln\left(\frac{1}{\ep}\right)+O(1).\label{hb6}
\end{eqnarray}
Combining \eqref{hb50} and \eqref{hb6}, we get  
\begin{eqnarray*}
	\frac{	\int_{\rr_+^N}\left| \Delta v_{\ep}\right| ^2\, dx}{\int_{\rr_+^N}\frac{v_{\ep}^2}{|x|^4}\, dx}=\frac{(N^2-4)^2}{16}+o(1) \bb{ as } \ep\to 0.
\end{eqnarray*}
 Since $v_\ep\in\H$, we get that $\gamma_H(\rn_+, \H)\leq \frac{(N^2-4)^2}{16}$. This completes the proof of Lemma \ref{lemma:extraonhalfspace}\qed\par 

\medskip\noindent To conclude this section, we discuss the model solutions for the homogeneous equation.
\begin{lem}\label{lemma:sol} Given $\alpha\in\rr$, for $x\in \rn_+$, $N\geq 5$, we define $v_\alpha(x):=x_1|x|^{-\alpha}$. 
			Then for $-N^2\leq \gg <\gg_{H,+}$, we have that
			\begin{equation}\label{eq:doublehardybi0}
				\left\{\begin{array}{ll}
					\Delta^2v_\alpha -\frac{\gg}{|x|^4} v_\alpha=0 &\hbox{ in } \rr_+^N, \\
					v_\alpha=\Delta v_\alpha = 0  &\hbox{ on } \partial \rr_+^N,
				\end{array}\right.
			\end{equation}
			if and only if $\alpha\in \{\alpha_{-}(\gg),\alpha_{+}(\gg),\beta_{-}(\gg),\beta_{+}(\gg)\}$ where
			\begin{equation}\label{p00}
				\alpha_{\pm}(\gg):=\frac{N-2}{2}\pm \frac{1}{2} \sqrt{N^2+4- 4\sqrt{N^2+\gg}},	
			\end{equation}
		and, 
		\begin{equation}\label{p01}
			\beta_{\pm}(\gg):=\frac{N-2}{2}\pm \frac{1}{2} \sqrt{N^2+4+ 4\sqrt{N^2+\gg}}.
		\end{equation}
		
		\end{lem}
\medskip\noindent \textit{Proof of Lemma \ref{lemma:sol}:} First, it follows from $-\Delta v=\alpha(N-\alpha)\frac{x_1|x|^{-\alpha}}{|x|^2} \bb{ on } \rr_+^N$ that 
			$$	\Delta^2v -\frac{\gg}{|x|^4} v=\left\{ \alpha(N-\alpha)(\alpha+2)(N-\alpha-2)-\gg\right\} \frac{v}{|x|^{4}}.$$
		To get our result, 	we want to solve the following equation:
			\begin{eqnarray}\label{p1}
				\alpha^4-2(N-2)\alpha^3+(N^2-6N+4)\alpha^2+(2N^2-4N)\alpha-\gg=0.
			\end{eqnarray}
			We denote $ a=1 \bb{ , } b=-2(N-2) \bb{ , } c=N^2-6N+4 \bb{ and } d=2N^2-4N$. Since $b^3-4abc+8a^2d=0$ then we can transform \eqref{p1} to biquadratic equation. We take $ \alpha:=t -\frac{b}{4a},$ and it follows from \eqref{p1} that 
			\begin{eqnarray}\label{p3}
				t^4-\left[ (N+2)^2+(N-2)^2 \right]\frac{t^2}{4}+\frac{(N+2)^2(N-2)^2}{16}-\gg=0.
			\end{eqnarray}
	It's easy to find the roots of \eqref{p3}. Since $ -N^2\leq \gg < \gg_{H,+}$, we find that 
			\begin{eqnarray*}	
				t_{\pm}( \gg):=\pm\frac{1}{2} \sqrt{N^2+4\pm 4\sqrt{N^2+\gg}}.
			\end{eqnarray*}
			Therefore, since $\alpha=t +\frac{N-2}{2}$ we get \eqref{p00}. Then, our fuction $v$ is a solution of equation \eqref{eq:doublehardybi0} when $\alpha\in\{ \alpha_{\pm}(\gg),\beta_{\pm}(\gg)\}$. This ends the proof of Lemma \ref{lemma:sol}. \qed\par 
	\medskip \noindent \textbf{Remarks about this lemma:} First, we remark that 
$$\beta_{-}(\gg)\leq \alpha_{-}(\gg)\leq \alpha_{+}(\gg)\leq \beta_{+}(\gg).$$
If $\gg=\gg_{H,+}$, then we have that 
$\alpha_{\pm}(\gg)=\frac{N-2}{2}.$ Hence, for $\gg\in[0,\gg_{H,+})$, we find that $$\alpha_{-}(\gg)\in \left[0,\frac{N-2}{2}\right)\bb{ and }\alpha_{+}(\gg)\in \left(\frac{N-2}{2}, N-2\right].$$
Moreover, note that if $\gg\in[0,\gg_{H,+})$, then   
$$\beta_{-}(\gg)\in\left(\beta_-(\gg_{H,+}),-2\right] \bb{ and }\beta_{+}(\gg)\in\left[N,\beta_+(\gg_{H,+})\right).$$
We also remark that $x_1|x|^{-\beta_{-}}$ is locally bounded  and $x_1|x|^{-\alpha_{-}}$ is the singular solution that is locally in $\H$. The following graph concerns the localizations of $\alpha_{\pm}$, $\beta_{\pm}$ when $\gg\in[0,\gg_{H,+})$:
\begin{figure}[!ht]
	\centering
	\includegraphics[scale=0.65, trim=40 20 40 20]{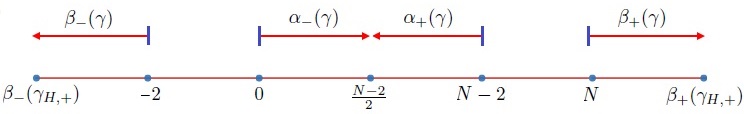}
\end{figure}

	\section{Existence of extremals for $Q_{\gg,s}(\rr_+^N)$: proof of Theorem \ref{theo1}}\label{sec:extr}
We fix $\gg < \gg_{H,+}=\frac{(N^2-4)^2}{16}$ and $0\leq s<4$. Recall that 
\begin{equation}\label{HShalfspace}
	Q_{\gg,s}(\rr_+^N)=\inf_{u\in \H\backslash \{0\}}\frac{\int_{\rr_+^N}\left( \left| \Delta u\right|^2-\gamma\frac{|u|^2}{|x|^4}\right) \, dx}{\left( \int_{\rr_+^N}\frac{|u|^{\crit}}{|x|^s}\,dx\right)^{\frac{2}{\crit}}}>0.
\end{equation} 
In order to prove the existence of extremals for \eqref{HShalfspace}, we proceed as in Ghoussoub-Robert \cite{GR} (see  Filippucci-Pucci-Robert \cite{FPR}):  these proofs rely on Lions's proof of the existence of extremals  for the Sobolev inequality \cite{Lio2}. We let $(\bar{u}_m)_{m\in \nn}\in \H$ be a minimizing sequence for $Q_{\gg,s}(\rr_+^N)$ in \eqref{HShalfspace} such that 
\begin{eqnarray*}
	\int_{\rr_+^N}\frac{|\bar{u}_m|^{\crit}}{|x|^s}\,dx =1 \bb{ and } \lim_{m\to+\infty}\left( \int_{\rr_+^N}\left| \Delta \bar{u}_m\right|^2\,dx-\gamma\int_{\rr_+^N}\frac{|\bar{u}_m|^2}{|x|^4}\,dx\right) = Q_{\gg,s}(\rr_+^N).
\end{eqnarray*}
For all $m\in \nn$, since $\int_{\rr_+^N}\frac{|\bar{u}_m|^{\crit}}{|x|^s}\,dx =1$, then, up to considering a subsequence,  there exists $\rho_m>0$ such that $\int_{\rr_+^N\cap B_{\rho_m}(0)}\frac{|\bar{u}_m|^{\crit}}{|x|^s}\,dx =\frac{1}{2}$. We define
$$ u_m(x):= \rho_m^{\frac{N-4}{2}}\bar{u}_m(\rho_m x) \bb{ for any } x\in \rr_+^N.$$
It is easy to check that $u_m \in \H$ for all $m\in \nn$, and 
\begin{eqnarray}\label{a00}
	\lim_{m\to +\infty}\left( \int_{\rr_+^N}\left| \Delta u_m\right|^2\,dx-\gamma\int_{\rr_+^N}\frac{|u_m|^2}{|x|^4}\,dx\right) =Q_{\gg,s}(\rr_+^N),
\end{eqnarray}
and 
\begin{eqnarray}\label{a0}
	\int_{\rr_+^N}\frac{|u_m|^{\crit}}{|x|^s}\,dx=1 \bb{ and } \int_{\rr_+^N\cap B_{1}(0)}\frac{|u_m|^{\crit}}{|x|^s}\,dx=\frac{1}{2}, 
\end{eqnarray}
for all $m\in \nn$. Since  $\gg <\gg_{H,+}$, there exists  $C>0$ such that 
\begin{equation}\label{eq:coercive}
	\int_{\rr_+^N}\left( \left| \Delta u\right|^2-\gamma\frac{|u|^2}{|x|^4}\right) \, dx\geq C\, \int_{\rr_+^N}\left| \Delta u\right|^2\, dx \bb{ for all } u\in \H.
\end{equation}
Therefore, with \eqref{a00}, there exists $C>0$ such that $C\,\vv u_m\vv^2  \leq Q_{\gg, s}(\rr_+^N)+o(1)$ as $m\to +\infty$. Hence, $(u_m)_{m\in \nn}$ is bounded in $\H$. As a consequence, up to the extraction of a subsequence, there exists $u\in \H$ such that 
\begin{eqnarray*}
	\left\{\begin{array}{ll}
		u_m \rightharpoonup  u \hbox{ weakly in } \H,\\
		u_m \to u \hbox{ strongly in } L^q_{loc}(\rr^N)\bb{ for any } 1\leq q <\crito:=\frac{2N}{N-4}.
	\end{array}\right.
\end{eqnarray*}
We define now the measures $ \nu_m, \la_m$ on $\rr^N$ as
\begin{eqnarray*}\label{eq:definitionmesures}
	\nu_m:=\frac{ |u_m|^{\crit}}{|x|^s}\mathbf{1}_{\rr_+^N}\, dx \hbox{ and }\la_m:= \left(\left| \Delta u_m\right|^2-\frac{\gg}{|x|^4}u_m^2 \right)\mathbf{1}_{\rr_+^N}\, dx.
\end{eqnarray*}
Using \eqref{a0} and \eqref{a00} yield 
\begin{eqnarray}\label{a000}
	\int_{\rr^N} d\nu_m =1 \bb{ and } \lim_{m\to +\infty} \int_{\rr^N}d\la_m=Q_{\gg,s}(\rr_+^N).
\end{eqnarray}
Up to extracting a subsequence, there exist two  measures $\la$, $\nu$ on $\rr^N$ such that
\begin{eqnarray}\label{eq:convmesures}
	  \la_m\rightharpoonup\la \bb{ and }   \nu_m \rightharpoonup \nu \bb{ weakly in the sens of measure as } m\to +\infty.
\end{eqnarray}
 We now apply Lions's first concentration-compactness Lemma \cite{Lio2} to the sequence of measures $(\nu_m)_{m\in \nn}$. Indeed, up to a subsequence, three situations may occur:
\begin{itemize}
	\item[(a)] (Compactness) There exists a sequence $(x_m)_{m\in \nn}$ in $\rr^N$ such that for any $\ep >0$ there exists $R_\ep>0$ with the property that 
	$$ \int_{B_{R_{\ep}}(x_m)} d\nu_m \geq 1- \ep \hbox{ for all } m\in \nn \bb{ large}.$$
	\item[(b)] (Vanishing) For all $R>0$ there holds
	$$ \lim_{m\to +\infty} \left(\sup_{x\in \rr^N} \int_{B_R(x)} \, d\nu_m \right) =0.$$
	\item[(c)] (Dichotomy) There exists $\alpha \in (0,1)$ such that for any $\ep >0$ there exists $R_\ep>0$ and a sequence $(x_m^\ep)_{m\in \nn}\in \rr^N$, with the following property: give $R^\prime >R_\ep$, there are non-negative measures $\nu_m^1$ and $\nu_m^2$ such that 
	\begin{align*}
		\hspace{1cm}	0\leq \nu_m^1&+\nu_m^2\leq \nu_m,\, Supp(\nu_m^1)\subset B_{R_\ep}(x_m^\ep), \, Supp(\nu_m^2)\subset \rr^N\backslash B_{R^\prime}(x_m^\ep),\\
		&	\nu_m^1=\nu_m \left|_{B_{R_\ep}(x_m^\ep)}\right., \hspace{0.8cm} \nu_m^2=\nu_m \left|_{\rr^N\backslash B_{R^\prime}(x_m^\ep)},\right.\\
		\lim_{m\to +\infty}&\sup \left(\left|\alpha-\int_{\rr^n} \, d\nu_m^1 \right| + \left|(1-\alpha)-\int_{\rr^n} \, d\nu_m^2 \right| \right)\leq \ep.
	\end{align*}
\end{itemize}
\begin{step}\label{step:lemma1}
	We claim that point $(a)$ (Compactness)  holds. In particular, we have that $\int_{\rr^N} \, d\nu =1$.
\end{step}
\smallskip\noindent \textit{Proof of Step \ref{step:lemma1}}:	Indeed, it follows from  \eqref{a0} that  point (b), does not hold. Assume by contradiction that  point (c) holds, that there exists $\alpha\in (0,1)$ such that (c) above holds. Taking $\ep=(m+1)^{-1}$, we can assume that, up to a subsequence, there exist $(R_m)_{m\in \nn}$ in $\rr_+$, $(x_m)_{m\in \nn}$ in $\rr^N$ and two sequence of non-negative measures, $(\nu_m^1)_{m\in \nn}$ and $(\nu_m^2)_{m\in \nn}$ such that $\lim\limits_{m\to +\infty} R_m=+\infty$ and
\begin{equation}\label{eq:convergencemesure}
	\left\{\begin{array}{ll}
		0\leq \nu_m^1+\nu_m^2\leq \nu_m&,\hspace{0.4cm} Supp(\nu_m^1)\subset B_{R_m}(x_m), \, Supp(\nu_m^2)\subset \rr^N\backslash B_{2R_m}(x_m),\\
			\nu_m^1=\nu_m \left|_{B_{R_m}(x_m)}\right.&, \hspace{0.4cm}  \nu_m^2=\nu_m \left|_{\rr^N\backslash B_{2R_m}(x_m)},\right.\\
		\lim\limits_{m\to +\infty}	\int_{\rr^N} \, d\nu_m^1=\alpha &,\hspace{0.12cm}\lim\limits_{m\to +\infty}	\int_{\rr^N} \, d\nu_m^2=1-\alpha.
	\end{array}\right.
\end{equation}
It follows from  \eqref{a000} and  \eqref{eq:convergencemesure} that
\begin{equation}\label{eq:convergenceDm}
	\lim_{m\to +\infty} \int_{B_{2R_m}(x_m)\backslash B_{R_m}(x_m)} \, d\nu_m=0.
\end{equation} 
Taking $\eta\in C_c^\infty(\rr^N)$ be such that $\eta_{|B_1(0)}\equiv1$, $\eta_{|B_2(0)^c}\equiv0$ and $0\leq \eta\leq 1$. For $m\in \nn$,  we define $\eta_m(x):=\eta\left( R_m^{-1}(x-x_m)\right)  \bb{ for all } x\in\rr^N.$ Using \eqref{eq:convergencemesure} and \eqref{eq:convergenceDm} 
\begin{eqnarray}
	1&=&\left( \int_{\rr^N} \eta_m^{\crit}\, d\nu_m^1+ \int_{\rr^N} (1-\eta_m)^{\crit}\, d\nu_m^2\right)^{\frac{2}{\crit}}+o(1)\nonumber\\
	&\leq& \left( \int_{\rr^N} \eta_m^{\crit}\, d\nu_m^1\right)^{\frac{2}{\crit}}+ \left( \int_{\rr^N} (1-\eta_m)^{\crit}\, d\nu_m^2\right)^{\frac{2}{\crit}}+o(1)\label{eq:etam1}.
\end{eqnarray}
On the other hand, it follows from \eqref{eq:convergencemesure} and \eqref{HShalfspace}  that 
\begin{eqnarray}
	&&\left( \int_{\rr^N} \eta_m^{\crit}\,d\nu_m^1\right)^{\frac{2}{\crit}}+ \left( \int_{\rr^N} (1-\eta_m)^{\crit}\,d\nu_m^2\right)^{\frac{2}{\crit}}+o(1)\nonumber\\
	&\leq& \left( \int_{\rr_+^N} \frac{|\eta_mu_m|^{\crit}}{|x|^s}\, dx\right)^{\frac{2}{\crit}}+ \left( \int_{\rr_+^N} \frac{|(1-\eta_m)u_m|^{\crit}}{|x|^s}\, dx\right)^{\frac{2}{\crit}}+o(1)\nonumber\\
	&\leq& Q_{\gg,s}(\rr_+^N)^{-1} \int_{\rr_+^N}\left( \left|\Delta (\eta_mu_m)\right|^2-\gamma\frac{|\eta_mu_m|^2}{|x|^4}\right) \,dx \nonumber\\
	&&+\,Q_{\gg,s}(\rr_+^N)^{-1}\int_{\rr_+^N}\left( \left| \Delta ((1-\eta_m)u_m)\right|^2-\gamma\int_{\rr_+^N}\frac{|(1-\eta_m)u_m|^2}{|x|^4}\right) \,dx+o(1).\label{eq:ineqdelta12}
\end{eqnarray}
\medskip\noindent{\bf Step \ref{step:lemma1}.1:} We claim that, as $m\to +\infty$, 
\begin{eqnarray}
	\int_{\rr_+^N}\left| \Delta (\eta_mu_m)\right|^2\,dx&=&	\int_{\rr_+^N}|\eta_m|^2\left|\Delta u_m\right|^2\,dx+o(1),\label{eq:delta1}\\
		\int_{\rr_+^N}\left| \Delta ((1-\eta_m)u_m)\right| ^2\,dx&=&\int_{\rr_+^N}|(1-\eta_m)|^2\left| \Delta u_m\right|^2\,dx+o(1).\label{eq:delta2}
\end{eqnarray}
\smallskip\noindent \textit{Proof of the claim:} We write for convenience $A_{+,m}:=\rr_+^N\cap \left( B_{2R_m}(x_m)\backslash B_{R_m}(x_m)\right) $. Since  $(u_m)_{m\in \nn}$ is bounded in $\H$, similarly to \eqref{eq:variationDelta} and \eqref{eq:Rep}, we get
\begin{eqnarray}\label{eq:variationDelta2}
	&&\int_{\rr_+^N}|\Delta (\eta_mu_m)|^2\,dx=\int_{\rr_+^N}|\eta_m|^2|\Delta u_m|^2 \,dx+O\left( \int_{A_{+,m}}\frac{|u_m|^2}{|x-x_m|^4}\,dx\right.\\
	&&+\left. \left( \int_{A_{+,m}}\frac{|u_m|^2}{|x-x_m|^4}\,dx\right)^{\frac{1}{2}}+\left( \int_{A_{+,m}}\frac{|u_m|^2}{|x-x_m|^4}\,dx\right)^{\frac{1}{4}}\right).\nonumber
\end{eqnarray}
We claim that,
\begin{equation}\label{eq:convergenceumx4}
	\lim_{m\to +\infty}\int_{A_{+,m}}\frac{u_m^2}{|x-x_m|^4}\,dx=0.
\end{equation}
\smallskip\textit{Proof.} 
Indeed, by Hölder's inequality and since $ N(1-\frac{2}{\crit})+\frac{2s}{\crit}=4$, we get that 
\begin{align}
	\int_{A_{+,m}}\frac{u_m^2}{|x-x_m|^4}\,dx&\leq R_m^{-4} \left[\int_{A_{+,m}}\, dx\right]^{1-\frac{2}{\crit}} \left[\int_{A_{+,m}} u_m^{\crit}\,dx\right]^{\frac{2}{\crit}}\nonumber\\
	&\leq c\,w_{N-1}^{1-\frac{2}{\crit}}{R_m}^{\frac{2s}{\crit}-4}\left[ \int_{R_m}^{2R_m}r^{N-1}\,dr\right]^{1-\frac{2}{\crit}}\left[\int_{A_{+,m}} \frac{u_m^{\crit}}{|x|^s}\,dx\right]^{\frac{2}{\crit}}\nonumber\\
	&\leq c\,w_{N-1}^{1-\frac{2}{\crit}}\,\left[ \int_{B_{2R_m}(x_m)\backslash B_{R_m}(x_m)} \, d\nu_m\right]^{\frac{2}{\crit}},\label{e16}
\end{align}
where $w_{N-1}$ is the volume of the canonical $(N-1)$--sphere. 
 Therefore, it follows from \eqref{eq:convergenceDm} that \eqref{eq:convergenceumx4}. This ends the proof of this claim.\qed\par

\smallskip\noindent Combining \eqref{eq:variationDelta2} and \eqref{eq:convergenceumx4}, we get the result. Similarly we  prove  \eqref{eq:delta2}. This ends the proof of  Step \ref{step:lemma1}.1.\qed\par

\medskip\noindent It follows from \eqref{eq:ineqdelta12}, \eqref{eq:delta1}, \eqref{eq:delta2} and \eqref{a000} that
\begin{eqnarray*}
	&&\left( \int_{\rr^N} \eta_m^{\crit}d\nu_m^1\right)^{\frac{2}{\crit}}+ \left( \int_{\rr^N} (1-\eta_m)^{\crit}d\nu_m^2\right)^{\frac{2}{\crit}}\\
	&\leq& Q_{\gg,s}(\rr_+^N)^{-1} \int_{\rr_+^N} \eta_m^2\left(|\Delta u_m|^2 -\gamma\frac{|u_m|^2}{|x|^4}\right)\, dx
\end{eqnarray*}
\begin{eqnarray*}
	&&+\,Q_{\gg,s}(\rr_+^N)^{-1} \int_{\rr_+^N} (1-\eta_m)^2\left(|\Delta u_m|^2 -\gamma\frac{|u_m|^2}{|x|^4}\right)\, dx+o(1)\\
	&\leq& Q_{\gg,s}(\rr_+^N)^{-1} \int_{\rr^N}\left(\eta_m^2+ (1-\eta_m)^2\right)\,d\la_m+o(1)\\ 
	&\leq& 1+2\, Q_{\gg,s}(\rr_+^N)^{-1}\int_{\rr_+^N} \eta_m(1-\eta_m) \frac{|u_m|^2}{|x|^4}\, dx+o(1)\\
	&\leq&1+O\left(\int_{A_{+,m} }\frac{|u_m|^2}{|x-x_m|^4}\, dx \right) +o(1).
\end{eqnarray*}
Letting $m \to +\infty$ and using \eqref{eq:convergencemesure}, \eqref{eq:etam1}  and \eqref{eq:convergenceumx4} yields $\alpha^{\frac{2}{\crit}}+(1-\alpha)^{\frac{2}{\crit}}=1$. This is impossible when $\alpha\in (0,1)$ and $\crit>2$. This contradiction proves Step \ref{step:lemma1}. \qed
\begin{step}\label{step:lemma2}
	There exists $I\subset \nn$ at most countable, and a family $\{x_i\}_{i\in I}\in \rr^N$ such that  
	\begin{equation}\label{eq:limitedenu}
		\nu=\frac{ |u|^{\crit}}{|x|^s}\mathbf{1}_{\rr_+^N}\, dx+\sum_{i\in I}\nu_i \delta_{x_i},
	\end{equation}
where $\nu_i:=\nu(x_i)>0$ for all $i\in I$. In particular, $\{x_i, i\in I\}\subset \{0\}$ when $s>0$. Moreover, there exists a bounded non–negative measure $\la_0\geq 0$ with no atoms (that is $\la_0(\{x\})=0$ for all $x\in \rr^N$) and 
\begin{equation}\label{eq:ll}
	\la=\la_0+\left( \left|\Delta u\right|^2- \gg \frac{ |u|^{2}}{|x|^4}\right) \mathbf{1}_{\rr_+^N}\, dx+\sum_{i\in I} \la_i \delta_{x_i}, 
\end{equation}
with  $\la_i=\la(\{x_i\})>0$ and  $\la_i \geq Q_{\gg,s}(\rr_+^N) \nu_i^{\frac{2}{\crit}}$.
\end{step}
\medskip\noindent \textit{Proof of Step \ref{step:lemma2}:} For $s=0$, \eqref{eq:limitedenu} is a consequence Lions's second concentration–compactness Lemma \cite{Lio2}. Take now $s>0$ so that $\crit<\frac{2N}{N-4}$, then $u_m\to u$ strongly in $L^{\crit}_{loc}(\rr_+^N)$. Therefore, we obtain that 
\begin{equation}\label{eq:ucrit}
	\nu= \frac{ |u|^{\crit}}{|x|^s}\mathbf{1}_{\rr_+^N}\, dx+\nu(\{0\})\delta_0.
\end{equation}
This proves \eqref{eq:limitedenu} in the case $s\geq 0$.\par
\smallskip\noindent We now prove \eqref{eq:ll} of Step \ref{step:lemma2}. We start by the following claim.\par 
	\medskip\noindent \noindent{\bf Step \ref{step:lemma2}.1:} We claim that 
	\begin{equation}\label{eq:controlenulambda}
		\left( \nu(\{x\})\right)^{\frac{2}{\crit}}\leq Q_{\gg,s}(\rr_+^N)^{-1}\la(x) \bb{ for all } x\in \rr^N.
	\end{equation}

\medskip\noindent \textit{Proof of the claim:} Indeed, $\varphi\in C^\infty(\rr^n)$ be such that $\varphi(x)=1$ for $x\in B_{1}(0)$, $\varphi(x)=0$ for $x\in \rr^N\backslash B_{2}(0)$ and $0\leq \varphi\leq 1$. Given $y\in\rr^N$ and $\delta >0$, we define $\varphi_\delta(x)=\varphi(\frac{x-y}{\delta})\bb{ for all } x\in \rr^N.$ Since $\varphi_\delta u_m\in \H$, the definition \eqref{eq:qgammas} yields
\begin{equation}\label{eq:Delta123}
	\left(\int_{\rr_+^N }\frac{|\varphi_\delta u_m|^{\crit}}{|x|^s}\,dx \right)^{\frac{2}{\crit}} \leq Q_{\gg,s}(\rr_+^N)^{-1}\int_{\rr_+^N}\left( \left| \Delta( \varphi_\delta u_m)\right| ^2-\gamma\frac{|\varphi_\delta u_m|^2}{|x|^4}\right) \,dx.
\end{equation}
As in the last proof of Step \ref{step:lemma1}.1 (see \eqref{eq:variationDelta2}), we have that 
\begin{eqnarray*}
	\int_{\rr_+^N}|\Delta (\varphi_\delta u_m)|^2\,dx=	\int_{\rr_+^N}|\varphi_\delta|^2|\Delta u_m|^2\,dx+R_{m,\delta}+o(1),
\end{eqnarray*}
where $o(1)\to 0$ as $m\to +\infty$ and,
\begin{eqnarray*}
	R_{m,\delta}:=O\left(\int_{A_{+,\delta}} \frac{|u_m|^{2}}{|x-y|^4}\,dx+\left( \int_{A_{+,\delta}} \frac{|u_m|^{2}}{|x-y|^4}\,dx\right)^\frac{1}{2}\right. \left.+ \left( \int_{A_{+,\delta}}\frac{|u_m|^2}{|x-y|^4}\,dx\right)^{\frac{1}{4}} \right),
\end{eqnarray*} 
where $A_{+,\delta}:=\rr_+^N\cap \left(B_{2\delta}(y)\backslash B_{\delta}(y) \right)$.
Therefore,  for all $\delta >0$, using \eqref{eq:Delta123} yields
\begin{align*}
\left(\int_{\rr^N }\left| \varphi_\delta\right| ^{\crit} \, d\nu_m \right)^{\frac{2}{\crit}}& \leq Q_{\gg,s}(\rr_+^N)^{-1}\int_{\rr^N}\varphi_\delta^2 \,d\la_m+R_{m,\delta}+o(1),
\end{align*}
letting $m\to +\infty$ and then $\delta \to 0$, we get that \eqref{eq:controlenulambda}. \qed \par
\smallskip\noindent Up to extraction,  let $\la^\prime$ be the weak limit of $|\Delta u_m|^2\mathbf{1}_{\rr_+^N}\, dx$ as $m\to +\infty$. Since $u_m\rightharpoonup u$ weakly in $\H$ as $m\to +\infty$, we get that $\la^\prime \geq |\Delta u|^2\mathbf{1}_{\rr_+^N}\, dx.$
Hence
\begin{equation}\label{eq:llprime}
	\la^\prime= \la_0 +\left| \Delta u\right|^2\mathbf{1}_{\rr_+^N}\, dx+\sum_{j\in J} \la^\prime(\{z_j\})\delta_{z_j}+\la^\prime_0\delta_0,
\end{equation}
where $\la_0\geq 0$ with no atoms, $z_j^{,}s$, $j\in J$ countable, and are the atoms of $\la^\prime$.
\noindent As above \eqref{eq:ucrit}, we have that there exists $L\geq 0$ such that 
\begin{equation}\label{eq:limiteux4}
	\frac{ |u_m|^{2}}{|x|^4}\mathbf{1}_{\rr_+^N}\, dx\rightharpoonup \frac{ |u|^{2}}{|x|^4}\mathbf{1}_{\rr_+^N}\, dx+L \, \delta_0,
\end{equation}
It follows from \eqref{eq:llprime} and \eqref{eq:limiteux4} that 
\begin{equation}\label{mino:ll}
	\la =\la_0+ \left( |\Delta u|^2- \gg \frac{ |u|^{2}}{|x|^4}\right) \mathbf{1}_{\rr_+^N}\, dx-\gg L\delta_0+\sum_{j\in J} \la^\prime(\{z_j\})\delta_{z_j}.
\end{equation}
First, using \eqref{mino:ll} and \eqref{eq:controlenulambda} yield
\begin{equation}\label{eq:nu0}
0<	\left( \nu(\{0\})\right)^{\frac{2}{\crit}} Q_{\gg,s}(\rr_+^N)\leq \la(\{0\}) =(\la^\prime_0(\{0\})-\gg L ).
\end{equation}
On the other hand, for $x_j\neq 0$ we have 
\begin{equation}\label{eq:nuj}
0<	\left( \nu(\{x_j\})\right)^{\frac{2}{\crit}} Q_{\gg,s}(\rr_+^N)\leq \la(\{x_j\}) =\la^\prime(\{x_j\}).
\end{equation}
From \eqref{mino:ll}, \eqref{eq:nu0} and \eqref{eq:nuj}, we obtain the result \eqref{eq:ll}. This proves Step \ref{step:lemma2}.\qed
\begin{step}\label{step:lemma12}
	We claim that one and only one of the two following situations occur:
	\begin{eqnarray*}
		&\hbox{either} 	\left\lbrace \nu =\frac{ |u|^{\crit}}{|x|^s}\mathbf{1}_{\rr_+^N}\, dx \hbox{ and } \int_{\rr_{+}^n} \frac{ |u|^{\crit}}{|x|^s} \, dx=1\right\rbrace \\
		& or\, \Big\{	\hbox{there exists $x_0\in \rr^N$ such that } \nu =\delta_{x_0} \hbox{ and } u\equiv0 \Big\} .
	\end{eqnarray*}

\end{step}
\medskip\noindent \textit{Proof of Step \ref{step:lemma12}:} Indeed, it follows from Step \ref{step:lemma1} that, 
\begin{align}\label{eq:i1}
	1&=\left( \int_{\rr^N} d\nu \right)^{\frac{2}{\crit}}= \left(\int_{\rr_+^N}\frac{|u|^{\crit}}{|x|^s}\, dx +\sum_{i\in I} \nu^i \int_{\rr^N} \delta_{x_i}\, dx  \right)^{\frac{2}{\crit}} \bb{ from \eqref{eq:limitedenu}}\nonumber\\
	&=\left(\int_{\rr_+^N}\frac{|u|^{\crit}}{|x|^s}\, dx +\sum_{i\in I} \nu^i   \right)^{\frac{2}{\crit}}\leq \left(\int_{\rr_+^N}\frac{|u|^{\crit}}{|x|^s}\, dx\right)^{\frac{2}{\crit}}+\sum_{i\in I}\nu_i^{\frac{2}{\crit}}. 
\end{align}
Now, using again \eqref{HShalfspace} and  \eqref{eq:controlenulambda} yield,
\begin{align}\label{eq:i2}
	&\left(\int_{\rr_+^N}\frac{|u|^{\crit}}{|x|^s}\, dx\right)^{\frac{2}{\crit}}+\sum_{i\in I}\nu_i^{\frac{2}{\crit}}\nonumber\\
	&\leq Q_{\gg,s}(\rr_+^N)^{-1}\left(\int_{\rr_+^N} \left( |\Delta u|^2-\gg \frac{u^2}{|x|^4}\right) \, dx+\sum_{i\in I} \la^i\right) \leq Q_{\gg,s}(\rr_+^N)^{-1} \int_{\rr^N} d\la
\end{align}
from \eqref{eq:ll}. Combining \eqref{eq:i1} and \eqref{eq:i2}, we have that $ \int_{\rr^N} d\la\geq Q_{\gg,s}(\rr_+^N)$.\par 
\smallskip \noindent We claim now that $ \int_{\rr^N} d\la\leq Q_{\gg,s}(\rr_+^N).$ Indeed, we let $f\in C^\infty(\rr^N)$ be such that $f(x)=0$ for $x\in B_{1}(0)$, $f(x)=1$ for $x\in \rr^N\backslash B_{2}(0)$ and $0\leq f\leq 1$. Given $\rho>0$, we let $f_\rho(x)=f(\rho^{-1}x)$ for all $x\in \rr^N$. So $(1-f_\rho^2)u_m\in  \H$  and therefore 
\begin{align*}
	\int_{\rr^N} (1-f_\rho^2)\, d\la_m&= \int_{\rr^N}  \, d\la_m- \int_{\rr_+^N} \left( f_\rho^2|\Delta u_m|^2-\gg \frac{|f_\rho u_m|^2}{|x|^4}\right) \, dx\\
	&=\int_{\rr^N}  \, d\la_m- \int_{\rr_+^N} \left( |\Delta (f_\rho u_m)|^2-\gg \frac{|f_\rho u_m|^2}{|x|^4}\right) \, dx\\
	&\qquad +\int_{\rr_+^N} \left( |\Delta (f_\rho u_m)|^2-f_\rho^2|\Delta u_m|^2\right) \, dx\\
	&\leq  \int_{\rr^N}  \, d\la_m +\int_{\rr_+^N} \left( |\Delta (f_\rho u_m)|^2-f_\rho^2|\Delta u_m|^2\right) \, dx\bb{ from \eqref{eq:coercive}}\\
	&\leq  Q_{\gg,s}(\rr_+^N)+R_{m,\rho}+o(1)\bb{ from \eqref{eq:variationDelta2}},
\end{align*}
where $o(1)\to 0 \bb{ as }m\to +\infty$, and
\begin{eqnarray*}
	R_{m,\rho}:=O\left(\int_{A_{+,\rho}} \frac{u_m^{2}}{|x-y|^4}\,dx+\left( \int_{A_{+,\rho}} \frac{u_m^{2}}{|x-y|^4}\,dx\right)^\frac{1}{2}\right. \left.+ \left( \int_{A_{+,\rho}}\frac{|u_m|^2}{|x-y|^4}\,dx\right)^{\frac{1}{4}} \right),
\end{eqnarray*} 
where $A_{+,\rho}:=\rr_+^N\cap\left( B_{2\rho}(0)\backslash B_{\rho}(0) \right)$. Therefore, letting $m\to +\infty$, and then $\rho\to +\infty$, and we then get this claim. 
\noindent This implies that $\int_{\rr^N}  \, d\la=Q_{\gg,s}(\rr_+^N)$. Therefore, it follows from \eqref{eq:i1} and \eqref{eq:i2} that 
$\left(\int_{\rr_+^N}\frac{|u|^{\crit}}{|x|^s}\, dx\right)^{\frac{2}{\crit}}+\sum_{i\in I}\nu_i^{\frac{2}{\crit}}= 1.$ By convexity, we have that one and only one term in \eqref{eq:limitedenu} is nonzero, then there exist $i_0\in I$ such that $x_0:=x_{i_0}$ and
$$\left\lbrace\nu^{i_0}=1 \hbox{ and } \int_{\rr_+^N}\frac{|u|^{\crit}}{|x|^s}\, dx=0\right\rbrace  \hbox{ or } \left\lbrace \nu^{i_0}=0  \hbox{ and } \int_{\rr_+^N}\frac{|u|^{\crit}}{|x|^s}\, dx=1\right\rbrace ,$$
with the equation \eqref{eq:limitedenu}, the exist $x_0\in \rr^N$ such that we get the claim of Step \ref{step:lemma12}.\qed\par
\begin{step}\label{step:extunozero}
	Suppose that $u\not\equiv 0$. We claim that that $u$ is an extremal for $Q_{\gg,s}(\rr_+^N)$.
\end{step}
\medskip\noindent \textit{Proof of Step  \ref{step:extunozero}:} Since $u\not\equiv0$, it follows from the previous Step that we have $\left(\int_{\rr_+^N}\frac{|u|^{\crit}}{|x|^s}\, dx\right)^{\frac{2}{\crit}}=1$. Using again the Hardy-Sobolev inequality \eqref{HShalfspace} yields,
\begin{align*}
	Q_{\gg,s}(\rr_+^N)\leq \int_{\rr_+^N}\left( |\Delta u|^2-\gg\frac{u}{|x|^4}\right)\, dx. 
\end{align*}
On the other hand, we have $u_m\rightharpoonup u$ as $m\to +\infty$ and we get that 
$$\int_{\rr_+^N}\left( |\Delta u|^2-\gg\frac{u^2}{|x|^4}\right)\, dx \leq \lim_{m\to +\infty}\inf\int_{\rr_+^N}\left( |\Delta u_m|^2-\gg\frac{u_m^2}{|x|^4}\right)\, dx=Q_{\gg,s}(\rr_+^N).$$
Therefore, we get the equality $I_{\gamma, s}^{\e}(u)=Q_{\gg,s}(\rr_+^N)$. That is  $u$ is an extremal for $Q_{\gg,s}(\rr_+^N)$. We obtain the result of  Step \ref{step:extunozero}.\qed 
\begin{step}\label{step:melsun}
	We suppose that $u\equiv 0$. Then, we have\begin{align*}
		s=0, \lim_{m\to +\infty}\int_{\rr_+^N}\frac{ |u_m|^2}{|x|^4}\, dx=0 \bb{ and } |\Delta u_m|^2\, dx \rightharpoonup Q_{\gg,0}(\rr_+^N) \delta_{x_0},
	\end{align*}
	as $m\to +\infty$ in the sense of measures.
\end{step}
\noindent\textit{Proof of step \ref{step:melsun}:} Since $u\equiv0$, and it follows from the Step \ref{step:lemma12} that  there exists $x_0\in \rr^N$ such that $\nu =\delta_{x_0}$.\par 
\smallskip\noindent We claim that $x_0\neq 0$. Indeed, if $x_0=0$, we get that $\int_{B_{1/2}(0)} \, d\nu=1$ which contradicts \eqref{a0}. Therefore  $x_0\neq 0$.  Since $u_m\rightharpoonup 0$ weakly in $\H$ as $m\to +\infty$, then for any $1\leq q< \frac{2N}{N-4}$, we have $u_m\to 0$ strongly in $L^q_{loc}(\rr_+^N)$.\par 
 \smallskip\noindent  We claim that $s=0$. 	Indeed, we argue by  contradiction and assume that $s>0$, then $\crit <\frac{2N}{N-4}$.  Let $r>0$, $\bb{since } x_0\neq 0 \bb{ and } u_m\to 0 \bb{ strongly in }  L^{\crit}_{loc}(\rr_+^N)$. Hence, we have $\lim\limits_{m\to +\infty} \int_{B_{r}(x_0)\cap\rr_+^N }\frac{|u_m|^{\crit}}{|x|^s}dx=0$,
 and, it follows from \eqref{eq:convmesures} and $\nu= \delta_{x_0}$ that $\lim\limits_{m\to +\infty} \int_{B_{r}(x_0)\cap\rr_+^N }\frac{|u_m|^{\crit}}{|x|^s}dx=1,$
 for all $r>0$ enough, a contradiction to our assumption. 
 
\smallskip\noindent Therefore $s=0$, and we prove the rest of this Step. Let $\rho>0$  and $f\in C^\infty(\rr^N)$ be such that $f(x)=0$ for $x\in B_{\rho}(x_0)$, $f(x)=1$ for $x\in \rr^N\backslash B_{2\rho}(x_0)$ and $0\leq f\leq 1$.
We now define, $\varphi:=1-f^2 \hbox{ and } \psi:=f\sqrt{2-f^2}.$
Clearly $\varphi,\psi\in C^\infty(\rr^N)$ and $\varphi^2+\psi^2=1$. It follows from \eqref{HShalfspace} and \eqref{eq:variationDelta2} that
\begin{eqnarray*}
	&&Q_{\gg,0}(\rr_+^N)\left(\int_{\rr_+^N }|\varphi u_m|^{2^{\star}_0}\,dx \right)^{\frac{2}{2^{\star}_0}}
	\leq \int_{\rr_+^N}\left( |\Delta( \varphi u_m)|^2\,dx-\gamma\int_{\rr_+^N}\frac{|\varphi u_m|^2}{|x|^4}\right) \,dx\\
	&&\leq\int_{\rr_+^N}\varphi^2\left( |\Delta u_m|^2\,dx-\gamma\int_{\rr_+^N}\frac{ |u_m|^2}{|x|^4}\right) \,dx+R_{m,\rho}+o(1),
\end{eqnarray*}
where $o(1)\to 0$ as $m\to +\infty$, and
\begin{eqnarray*}
	R_{m,\rho}:=O\left(\int_{A_{+,\rho}} \frac{|u_m|^{2}}{|x-x_0|^4}\,dx+\left( \int_{A_{+,\rho}} \frac{|u_m|^{2}}{|x-x_0|^4}\,dx\right)^\frac{1}{2}\right. \left.+ \left( \int_{A_{+,\rho}}\frac{|u_m|^2}{|x-x_0|^4}\,dx\right)^{\frac{1}{4}} \right),
\end{eqnarray*} 
where $A_{+,\rho}:=\rr_+^N\cap\left( B_{2\rho}(x_0)\backslash B_{\rho}(x_0) \right)$. Using $u_m\to 0$ in $L^2_{loc}(\rr^N)$ yields $R_{m,\rho}=o(1)$ as $m\to +\infty$. And, so by $\varphi^2=1-\psi^2$ 
\begin{align}
		Q_{\gg,0}(\rr_+^N)\left(\int_{\rr_+^N }|\varphi u_m|^{2^{\star}_0}\,dx \right)^{\frac{2}{2^{\star}_0}}&\leq \int_{\rr_+^N}\left( |\Delta u_m|^2\,dx-\gamma\int_{\rr_+^N}\frac{ |u_m|^2}{|x|^4}\right) \,dx\nonumber\\
	&\quad-\int_{\rr_+^N}\psi^2\left( |\Delta u_m|^2\,dx-\gamma\int_{\rr_+^N}\frac{ |u_m|^2}{|x|^4}\right) \,dx+o(1), \label{eq:phipsi}
\end{align}
$\bb{ as }m\to +\infty.$ It follows from \eqref{eq:convmesures} and Step \ref{step:lemma12} that
\begin{equation}\label{eq:phium}
	\left(\int_{\rr_+^N }|\varphi u_m|^{2^{\star}_0}\,dx \right)^{\frac{2}{2^{\star}_0}}=\left( |\varphi (x_0)|^{\crito}+o(1)\right)^{\frac{2}{2^{\star}_0}}=1+o(1) \bb{ as } m\to +\infty.
\end{equation}
Plugging \eqref{eq:phium} into \eqref{eq:phipsi} we get that
\begin{align*}
	Q_{\gg,0}(\rr_+^N)+o(1)&\leq \int_{\rr_+^N}\left( |\Delta u_m|^2-\gamma\frac{ |u_m|^2}{|x|^4}\right) \,dx\nonumber\\
	&-\int_{\rr_+^N}\psi^2\left( |\Delta u_m|^2-\gamma\frac{ |u_m|^2}{|x|^4}\right) \,dx+o(1),\label{eq:new1}
\end{align*}
as $m\to +\infty$. From \eqref{a00}, we have 
\begin{equation}\label{eq:psi}
	\int_{\rr_+^N}\psi^2\left( |\Delta u_m|^2-\gamma\frac{ |u_m|^2}{|x|^4}\right) \,dx\leq o(1) \hbox{ as } m \to +\infty.
\end{equation}
As in the proof of \eqref{eq:delta1}, we obtain that $	\int_{\rr_+^N} |\Delta( \psi u_m)|^2\,dx=\int_{\rr_+^N}\psi^2 |\Delta u_m|^2\,dx+o(1)$ as $m\to +\infty$. Plugging this expansion into \eqref{eq:psi} yields
$$\lim\limits_{m\to+\infty}\int_{\rr_+^N}\left( |\Delta (\psi u_m)|^2-\gamma|x|^{-4} |\psi u_m|^2\right) \,dx=0.$$
 Hence, by the coercivity \eqref{eq:coercive}, we get  
 \begin{equation}\label{convergenceDelta}
	\lim_{m\to +\infty} \vv \Delta (\psi u_m) \vv_2=0.
\end{equation}
With the result of Lemma \ref{lemma:extraonhalfspace}, we have that $$	\lim\limits_{m\to +\infty}\int_{\rr_+^N}\frac{ |\psi u_m|^2}{|x|^4}\, dx=0.$$ We then have $\lim\limits_{m\to +\infty}\int_{\rr_+^N\backslash B_{2\rho}(x_0)}\frac{ |u_m|^2}{|x|^4}\, dx=0.$ Moreover, taking $\rho>0$ small enough and since $u_m\to 0$ in $L^2_{loc}(\rr_+^N)$ around $x_0\neq 0$, we can conclude that  $$\lim\limits_{m\to +\infty}\int_{\rr_+^N}\frac{ |u_m|^2}{|x|^4}\, dx=0,$$
which implies  by \eqref{a00} that $\lim\limits_{m\to +\infty}\vv \Delta u_m\vv^2=Q_{\gg,0}(\rr_{+}^N)$. Hence, using \eqref{convergenceDelta}  yields the third part of claim.  Step \ref{step:melsun} is proved.\qed
\begin{step}\label{step:mq00rn}
	We now claim that, if $u\equiv0$, then $s=0$ and $Q_{\gg,s}(\rr_+^N)=S_N$,
where $S_N$ is defined in \eqref{eq:bestconstantsoblevrn}.
\end{step}
\noindent\textit{Proof of step \ref{step:mq00rn}:} We have already seen that $s=0$. Since $u_m\in \H\subset \tilde{H}_0(\rn)$, we have that 
\begin{eqnarray*}
S_N\left(\int_{\rr^N} |u_m|^{2^{\star}_0} \, dx \right)
	&\leq&\int_{\rr^N} |\Delta u_m|^2\, dx\\
	&\leq&\int_{\rr_+^N} \left( |\Delta u_m|^2-\gg \frac{u_m^2}{|x|^4}\right) \, dx\\
	&&+\int_{\rr^N\backslash\rr_+^N} |\Delta u_m|^2\, dx+\,\gg\int_{\rr_+^N}\frac{u_m^2}{|x|^4} \, dx,
\end{eqnarray*}
 so with \eqref{a00} and by the result of Step \ref{step:melsun}, we get 
\begin{eqnarray*}
	S_N\left(\int_{\rr^N} |u_m|^{2^{\star}_0} \, dx \right)& \leq& Q_{\gg,0}(\rr_+^N) +o(1).
\end{eqnarray*}
It follows then from \eqref{a0} that $S_N\leq Q_{\gg,0}(\rr_+^N) +o(1)$. Letting $m\to +\infty$  we obtain that $S_N\leq Q_{\gg,0}(\rr_+^N)$. Conversely, it follows from the computations of Proposition \ref{prop:nonext} below that $Q_{\gg,0}(\rr_+^N)\leq S_N$. Hence, we have $Q_{\gg,0}(\rr_+^N)= S_N$. This proves Step \ref{step:mq00rn} \qed
\begin{step}\label{step:casegammaneagtive}
	We assume that $s=0$ and $\gg\leq 0$. Then, we have that $Q_{\gg,0}(\rr_+^N)= S_N$, where $S_N$ is defined in \eqref{eq:bestconstantsoblevrn}.
\end{step}
\noindent\textit{Proof of step \ref{step:casegammaneagtive}:} Indeed, since $\gg \leq 0$, it follows from Proposition \ref{prop:nonext} when $\Omega= \rr_+^N$ is smooth domain that  $Q_{\gg,0}(\rr_+^N)= S_N.$ \qed
\begin{step}\label{step:casegammapositive}
	Taking $\{s=0, \gg>0 \bb{ and } N\geq 8\}$, we claim that
	\begin{equation*}\label{eq:conditonexistence}
		Q_{\gg,0}(\rr_+^N)< S_N.
	\end{equation*}
Therefore, we get that $Q_{\gg,0}(\e)$ is attained.
\end{step}
\noindent\textit{Proof of step \ref{step:casegammapositive}:} We fix $x_0\in \rr_+^N$ such that $x_0\neq 0$. We define $U(x):=(1+|x|^2)^{-\frac{N-4}{2}}$ for all $x\in\rn$. It follows from Lions \cite{Lio1,Lio2} that $U\in \tilde{H}_0(\rn)=D^{2,2}(\rn)$ is an extremal for \eqref{eq:bestconstantsoblevrn}, that is
$$S_N:=\frac{\int_{\rn}\left| \Delta U\right| ^2\, dx}{\left(\int_{\rn}|U|^{\crito}\, dx \right)^{\frac{2}{\crito}}}.$$
Let $\eta \in C_c^{\infty}(\rr_+^N)$ and $0<\delta<|x_0|/2$ be such that $\eta(x)=1$ for $x\in B_\delta(x_0)$. We consider the test function  
\begin{equation}\label{def:Ue}
	U_\ep(x):=\eta(x)u_\ep(x) \hbox{ for all } x\in \rr_+^N\bb{ and } \ep >0,
\end{equation}
\begin{equation*}\label{eq:defuepsilon}
\hbox{where }	u_\ep(x):=\ep^{-\frac{N-4}{2}}U\left(\frac{x-x_0}{\ep}\right)=\left( \frac{\ep}{\ep^2+|x-x_0|^2}\right)^{\frac{N-4}{2}}\hbox{ for all }x\in\rn.
\end{equation*} 

\medskip\noindent{Step \ref{step:casegammapositive}.1:}
We claim as $\ep \to 0$ that:
\begin{eqnarray}\label{eq:i1ep}
	\int_{\rr_+^N} \frac{U_\ep^2}{|x|^4}\, dx=\left\{\begin{array}{ll}
		\frac{\ep^{4}}{2}|x_0|^{-4}\int_{\rr^N}U^2\, dx+o(\ep^{4}) &\bb{ if  } N\geq 9,\\
		w_{7}\,\ep^4	\ln(\frac{1}{\ep})|x_0|^{-4}+o\left(\ep^4	\ln\left(\frac{1}{\ep}\right)\right) &\bb{ if } N=8,\\
		O(\ep^{N-4})&\bb{ if } N=5,6,7,
	\end{array}\right.
\end{eqnarray}
where $w_7$ is the volume of the canonical $7$–sphere.\par
\medskip\noindent \textit{Proof of the claim:} Indeed,
for $\delta >0$ we begin by noticing that 
\begin{equation}\label{eq:Uepx4}
	\int_{\rr_+^N} \frac{U_\ep^2}{|x|^4}\, dx
	=|x_0|^{-4} I_{1,\ep}+I_{2,\ep}+O(\ep^{N-4}) \bb{ as } \ep \to 0,
\end{equation}
where:
\begin{eqnarray*}
	I_{1,\ep}:=\int_{B_\delta(x_0)\cap\e}u_\ep^2\, dx \bb{ and } I_{2,\ep}:=\int_{B_\delta(x_0)\cap\e}\left(\frac{1}{|x|^4} -\frac{1}{|x_0|^4}\right)  u_\ep^2\, dx.
\end{eqnarray*}
First, we claim that 
\begin{eqnarray}\label{eq:i11ep}
	I_{1,\ep}=\left\{\begin{array}{ll}
		\frac{\ep^{4}}{2}\int_{\rr^N}U^2\, dx+O(\ep^{N-4}) &\bb{ if  } N\geq 9,\\
		w_{7}\,\ep^4	\ln(\frac{1}{\ep})+O(\ep^4) &\bb{ if } N=8,\\
		O(\ep^{N-4})&\bb{ if } N=5,6,7.
	\end{array}\right.
\end{eqnarray}
\begin{proof} For $N\geq 9$, $U\in L^2(\rn)$, and we get as $\ep \to 0$ that
\begin{equation*}
		I_{1,\ep}=\frac{\ep^{4}}{2}\int_{\rr^N}U^2\, dx+O(\ep^{N-4}).
\end{equation*}
Take now the case $N=8$. It follows from the change of variable in polar coordinates that
\begin{align*}
	I_{1,\ep}	&=\ep^4\, w_7\left[O(1)+\int_{1}^{\delta\ep^{-1}}\frac{1}{r}\, dr+\int_{1}^{\delta\ep^{-1}}r^{7}\left[ \frac{1}{\left( 1+r^2\right)^4 }-\frac{1}{r^8}\right]\, dr \right]\\ 
	&=\ep^4w_7 \ln\left(\frac{\delta}{\ep}\right)+O(\ep^4).
\end{align*}
As one checks, we have that $I_{1,\ep}=O(\ep^{N-4})$ for $N=5,6,8$. This proves \eqref{eq:i11ep}.
\end{proof} 
\noindent Next, we claim as $\ep \to 0$ that 
\begin{eqnarray}\label{eq:i2ep}
	I_{2,\ep}=\left\{\begin{array}{ll}
		o(\ep^4) &\bb{ if  } N\geq 9,\\
		o(\ep^4\ln(\frac{1}{\ep})) &\bb{ if } N=8,\\
		O(\ep^{N-4})&\bb{ if } N=5,6,7.
	\end{array}\right.
\end{eqnarray}
\begin{proof} It follows from the definition of $u_\ep$ that
	\begin{align*}
		\left| I_{2,\ep}\right| \leq \ep^{N-4}\, \int_{B_\delta(0)\cap\e}\frac{|A(x)|}{\left( \ep^{2}+|x|^{2}\right)^{N-4} }\, dx,
	\end{align*}
	where $A(x):=\frac{1}{|x+x_0|^4}-\frac{1}{|x_0|^4}$.  Fix $\alpha \in (0,\delta)$, with a change of variables we write 
	\begin{equation}\label{eq:ineq1}
			\ep^{N-4}\, \int_{(B_{\delta}(0)\backslash B_\alpha(0))\cap\e}\frac{|A(x)|}{\left( \ep^{2}+|x|^{2}\right)^{N-4} }\, dx
		 \leq \ep^{N-4}\, \vv A\vv_{\infty} \frac{\delta^{N}w_{N-1}}{\alpha^{2(N-4)}}:=C_{1,\alpha}\ep^{N-4}.
	\end{equation}
On the other hand, we have that 
\begin{equation*}
	\ep^{N-4}\, \int_{ B_\alpha(0)\cap \e}\frac{|A(x)|}{\left( \ep^{2}+|x|^{2}\right)^{N-4} }\, dx\\
	\leq  \ep^{4}w_{N-1}\left(\sup_{|x|\leq \alpha}|A(x)|\right) \int_{ 0}^{\alpha	 \ep^{-1}}\frac{r^{N-1}}{\left( 1+r^{2}\right)^{N-4} }\, dr.
\end{equation*}
From here, we split the proof in three cases:\par
\smallskip \noindent \textbf{Case 1:} If  $N\geq 9$, we have $U\in L^2(\rr^N)$. In this case, we can write 
\begin{eqnarray}\label{eq:ineq2}
&&	\ep^{N-4}\, \int_{ B_\alpha(0)\cap\e}\frac{|A(x)|}{\left( \ep^{2}+|x|^{2}\right)^{N-4} }\, dx\nonumber\\
&&	\leq \ep^4\left(  \sup_{|x|\leq \alpha}|A(x)|\right) \left[ \int_{\rr^N}|U|^2\, dx -\frac{w_{n-1}\ep^{N-8}}{(N-8)\alpha^{N-8}}\right].
\end{eqnarray}
Combining \eqref{eq:ineq1} and \eqref{eq:ineq2}, then there exist $C_{2,\alpha}, C>0$ such that 
\begin{equation*}
	|I_{2,\ep}| \leq C_{2,\alpha}\ep^{N-4}+C\ep^4 \left( \sup_{|x|\leq \alpha}|A(x)|\right).
\end{equation*}
	Take now $\theta>0$. Then there exists $\alpha_0>0$ such that for all $\alpha<\alpha_0$ we have that $C\left(\sup_{|x|\leq \alpha}|A(x)|\right)  \leq  \frac{\theta}{2}$. On the other hand, we have $\lim\limits_{\ep\to 0}C_{2,\alpha} \ep^{N-8}=0$, then there exists $\ep_0:=\ep(\theta)>0$ such that $C_{2,\alpha}\ep^{N-8}\leq \frac{\theta}{2}$. Then, we have $I_{2,\ep}=o(\ep^4)$ as $\ep\to 0$. \par 
\smallskip \noindent \textbf{Case 2:} The proof of case $N=8$ is similar to the proof of \eqref{eq:i11ep}. \end{proof}  
\noindent Plugging \eqref{eq:i11ep} and  \eqref{eq:i2ep} into \eqref{eq:Uepx4}, we obtain that \eqref{eq:i1ep}. This proves Step \ref{step:casegammapositive}.1. \qed\par
\noindent For $N\geq 5$, it is also classical as $\ep \to 0$ that:
\begin{eqnarray}
	\int_{\rr_+^N} \left|\Delta U_\ep (x) \right|^2\, dx&=& \int_{\rr^N}\left| \Delta U\right|^2\, dx+O(\ep^{N-4}),\label{eq:la}\\	\int_{\rr^N_+}| U_{\ep}|^{\crito}\, dx&=& \int_{\rr^N}| U|^{\crito}\,  dx+O(\ep^{N}).\label{eq:U2**}
\end{eqnarray}
\noindent Combining \eqref{eq:i1ep}, \eqref{eq:la} and \eqref{eq:U2**}, we have that
\begin{eqnarray*}
		I_{\gamma, 0}^{\e}(U_\ep)=\left\{\begin{array}{ll}
		S_N-\gg|x_0|^{-4}c\,\ep^{4}+o(\ep^{N-4}) &\bb{ if  } N\geq 9,\\
		S_N-\gg|x_0|^{-4}	c\,\ep^4\ln(\frac{1}{\ep})+o(\ep^4	\ln(\frac{1}{\ep})) &\bb{ if } N=8,
	\end{array}\right.
\end{eqnarray*}
where $c$ is a positive constant. Since $\gg >0$ and $x_0\neq 0$, then $Q_{\gg,0}(\rr_+^N)< S_N$. Therefore, it follows from the Step \ref{theo1} that $u\not\equiv 0$, and we have $u$ is a extremal for $Q_{\gg,0}(\e)$. This ends Step \ref{step:casegammapositive}.\qed

\medskip\noindent All these cases end the proof of Theorem \ref{theo1}.

\section{Proof of Theorem \ref{theo2}}\label{sec:theo2}
In this section, we use the existence of extremals for Hardy-Sobolev inequality, established in Section \ref{sec:extr} to prove that there exists a nontrivial weak solution for double critical equation \eqref{eq:doublehardybi}.\par 
\smallskip \noindent For any functional $G\in C^1(X,\rr)$ where $(X,\vv \cdot\vv)$ is a Banach space, we say that $(u_m)_{m\in \nn}\in X$ is a Palais-Smale sequence of  $G$ if  there exists $\beta\in\rr$ such that
\begin{eqnarray*}
	G(u_m)\to \beta  \hbox{ and } G^{\prime}(u_m)\to 0 \hbox{ in } X^{\prime} \hbox{ as } m\to +\infty.
\end{eqnarray*}
Here, we say that the Palais-Smale sequence is at level $\beta$. The main tool is the Mountain-Pass lemma of Ambrosetti-Rabinowitz \cite{AR}:
\begin{theo}[Mountain-Pass lemma \cite{AR}]\label{thAR}
	We consider $G\in C^{1}(X,\rr)$ where $(X,\vv \cdot\vv)$ is a Banach space. We assume that $G(0)=0$ and that
	\begin{itemize}
		\item There exist $\lambda, r>0$ such that $G(u)\geq \lambda$ for all $u\in X$ such that $\vv u \vv=r$,
		\item There exists $u_0$ in $X$ such that $\lim sup_{t\to +\infty}G(tu_0)<0$.
	\end{itemize}
	We consider $t_0>0$ sufficiently large such that $\vv t_0u_0\vv >r$ and $G(t_0u_0)<0$, and
	$$\beta=\inf_{c\in \Gamma} \sup_{t\in [0,1]} G(c(t)),$$
	
	$$ \hbox{where }\Gamma:=\{c\in C^0([0,1],X) \hbox{ s.t. } c(0)=0, \, c(1)=t_0u_0\}.$$
	Then, there exists a Palais-Smale sequence at level $\beta$ for $G$. Moreover, we have that $ \beta \leq \sup_{t\geq 0} G(tu_0)$.
\end{theo}
\noindent We  define the energy functional noted by $\E $ 
\begin{equation*}
	\E (u):=\frac{1}{2}\int_{\e}\left( \left| \Delta  u\right| ^2-\gg\frac{u^2}{|x|^4}\right)\, dx-\frac{1}{\crit}\int_{\e}\frac{|u|^{\crit}}{ |x|^s}\,dx- \frac{1}{\crito}\int_{\e} |u|^{\crito}\, dx,
\end{equation*}
for any  $u\in \H$. Any weak solution to \eqref{eq:doublehardybi} is a critical point of $\E $. In the sequel, since $\gg <\gg_{H,+}$, then \eqref{eq:coercive} holds.
\begin{prop}\label{prop0}
	We assume that $\gamma<\gamma_{H,+}$. Fix $u_0\in \H$ such that $u_0\geq 0$, $u_0\not\equiv 0$. Then there exists a sequence $(u_m)_{m\in\nn}\in \H$ that is a Palais-Smale sequence for $\E $ at level $\beta$ such that $0<\beta\leq \sup_{t\geq 0}\E (tu_0)$.
\end{prop}
\noindent{\it Proof of Proposition \ref{prop0}:} Indeed, clearly $\E \in C^1(\H)$.  Note that $\E(0)=0$. It follows from \eqref{eq:coercive} and the Sobolev and Hardy-Sobolev embeddings that there exist $c_0,c_1,c_2>0$ such that
\begin{equation}\label{1}
	\E(u)\geq c_0 \vv u \vv^2 -c_1\vv u\vv^{\crit}-c_2\vv u\vv^{\crito} \hbox{ for all } u\in \H.
\end{equation}
Define $f(r)=r^2 \left[ c_0-c_1r^{\crit-2}-c_2r^{\crito-2}\right] :=r^2 g(r)$ and since $\crit,\crito>2$ we have $g(r)\to c_0$ as $r\to 0$. Then there exists $r_0>0$ such that $r<r_0$, we have $g(r)>\frac{c_0}{2}$. Therefore, for all $u\in\H$ such that $\vv u \vv =\frac{r_0}{2}$ and by \eqref{1}, we have $\E(u)\geq \frac{c_0r_0^2}{8}:=\lambda$. We fix $u_0\in\H$, $u_0\not\equiv 0$, and
\begin{eqnarray*}
	\E(tu_0)&=&\frac{t^2}{2}\int_{\e}\left( \left| \Delta u_0\right| ^2-\frac{\gg}{|x|^4}u_0^2\right) dx\\
	&-&\frac{t^{\crit}}{\crit}\int_{\e}\frac{|u_0|^{\crit}}{|x|^s}dx
	-	\frac{t^{\crito}}{\crito} \int_{\e}|u_0|^{\crito}\,  dx\\
	&:=&\frac{t^2}{2}R_1
	-\frac{t^{\crit}}{\crit}R_2
	-	\frac{t^{\crito}}{\crito}R_3\leq t^{\crit}\left( \frac{t^{2-\crit}}{2}R_1
	-R_2\right),
\end{eqnarray*}
where $R_1,R_2>0$ and $R_3\geq 0$. Since $\crit>2$, we have $\E (tu_0)\to -\infty$ as $t\to +\infty$. Then $\limsup_{t\to +\infty} \E (tu_0)<0$. 	We consider $t_0>0$ large such that $\vv t_0u_0\vv >r$ and $\E (t_0u_0)<0$. For $t\in[0,1]$, we have $\E (c(t))\geq \lambda$ and then there exists 
\begin{eqnarray*}
	\beta:=\inf_{c\in \Gamma} \sup \E (c(t))\geq \lambda >0.
\end{eqnarray*}
We apply Theorem \ref{thAR} to get the expected Palais-Smale sequence. This ends the proof of Proposition  \ref{prop0}.\qed
\begin{prop}\label{prop00}
Suppose $0<\gg <\gg_{H,+}$, $0\leq s<4$ and $N\geq 8$.	 Then there exists a sequence $(u_m)_{m\in\nn}\in \H$ that is a Palais-Smale sequence for $\E $ at level $\beta$ such that 
\begin{equation}\label{eq:levelbeta}
	0<\beta <\beta^{\star}:=\min\left\lbrace\frac{2}{N}\,Q_{\gg,0}(\e)^{\frac{N}{4}}, \frac{4-s}{2(N-s)}Q_{\gg,s}(\e)^{\frac{N-s}{4-s}}\right\rbrace .
\end{equation}

\end{prop} 
\medskip \noindent\noindent{\it Proof of Proposition \ref{prop00}:} From  Theorem \ref{theo1}, we know that there exists an  extremal $u_0\in\H$ for $Q_{\gg,0}(\e)$ whenever $\gg>0$ and $N\geq 8$. It follows then from Proposition \ref{prop0} that there exists a sequence $u_m\in \H$  a Palais-Smale sequence for $\E$ at level $\beta$ such that $$\beta\leq \sup_{t\geq 0}\E (tu_0)\leq \sup_{t\geq 0} f_1(t),$$
where: $$f_1(t):= \frac{t^2}{2}\int_{\e}\left( \left| \Delta  u_0\right| ^2-\frac{\gg}{|x|^4}u_0^2\right)\, dx- \frac{t^{\crito}}{\crito}\int_{\e} |u_0|^{\crito}\, dx \hbox{ for all } t>0.$$
Simple computations yield that $f_1(t)$ attains its maximum at the point $$t_{max}=\left(\frac{\int_{\e}\left( \left| \Delta  u_0\right| ^2-\frac{\gg}{|x|^4}u_0^2\right)\, dx}{\int_{\e} |u_0|^{\crito}\, dx} \right)^{\frac{1}{\crito-2}}. $$
Therefore, since $\crito=\frac{2N}{N-4}$ and $u_0$ is an extremal for $Q_{\gg,0}(\e)$, we have that 
\begin{align*}
\sup_{t\geq 0} f_1(t)&=\left[\frac{1}{2}-\frac{1}{\crito} \right]\left(\frac{\int_{\e}\left( \left| \Delta  u_0\right| ^2-\frac{\gg}{|x|^4}u_0^2\right)\, dx}{\left( \int_{\e} |u_0|^{\crito}\, dx\right)^\frac{2}{\crito} } \right)^{\frac{\crito}{\crito-2}}=\frac{2}{N}Q_{\gg,0}(\e)^{\frac{N}{4}}.
\end{align*}
Thus, $\beta \leq \sup_{t\geq 0} f_1(t)=\frac{2}{N}\,Q_{\gg,0}(\e)^{\frac{N}{4}}.$ We now prove that this inequality is strict. Assume by contradiction that
$$\sup_{t\geq 0}\E (tu_0)= \sup_{t\geq 0} f_1(t),$$
and we consider $t_1>0$ where $\sup_{t\geq 0}\E (tu_0)$ is attained. We obtain that
$$f_1(t_1)-\frac{t_1^{\crit}}{\crit}\int_{\e}\frac{|u_0|^{\crit}}{|x|^s}\, dx = f_1(t_{max}),$$
this give us $f_1(t_1)-f_1(t_{max})>0$ (because $t_1>0$). Contradiction with $t_{max}$ is a maximum point of $f_1(t)$. Therefore, we have $\beta <\frac{2}{N}\,Q_{\gg,0}(\e)^{\frac{N}{4}}$. Similar, we can get $\beta < \frac{4-s}{2(N-s)}Q_{\gg,s}(\e)^{\frac{N-s}{4-s}}$ whenever $s>0$. Thus, we can define $\beta^{\star}$ as in \eqref{eq:levelbeta} such that $0<\beta<\beta^{\star}$. This proves Proposition \ref{prop00}.\qed
\begin{prop}\label{prop2} 
Suppose $0<\gg < \gg_{H,+}$, $0\leq s<4$ and $N\geq 8$. We assume that $(u_m)_{m\in \nn}$ is a Palais-Smale sequence of $E$ at energy level $\beta \in (0, \beta^\star)$. If $u_m \rightharpoonup  0 \hbox{ weakly in } \H$ as $m\to +\infty$, then there exists $\ep:=\ep(N,\gamma,s,\beta)>0$  such that 
\begin{eqnarray*}
\hbox{ either }\lim_{m\to +\infty}\sup \int_{B_r(0)} \left| u_m\right|^{\crito}\, dx  =\lim_{m\to +\infty}\sup \int_{B_r(0)} \frac{\left| u_m\right|^{\crit}}{|x|^s}\, dx=0; 
\end{eqnarray*}
\begin{eqnarray*}
\hbox{ or }	\lim_{m\to +\infty}\sup \int_{B_r(0)} \left| u_m\right|^{\crito}\, dx\geq\ep  \hbox{ and } \lim_{m\to +\infty}\sup \int_{B_r(0)} \frac{\left| u_m\right|^{\crit}}{|x|^s}\, dx\geq \ep,
\end{eqnarray*}
for every $r>0$.
\end{prop}

\noindent{\it Proof of Proposition \ref{prop2}:} Indeed, the proof of this proposition is divided into several steps:

\begin{step}\label{step:00}
For $s\in(0,4)$. Let $K$ be an arbitrary compact set in $\e\backslash \{0\}$, we claim 
	\begin{align*}
		&\lim_{m\to +\infty} \int_{K}\frac{u_m^{\crit}}{|x|^s}\, dx=\lim_{m\to +\infty}\int_{K}\frac{u_m^2}{|x|^4}\, dx=0,\\
		&\lim_{m\to +\infty}\int_K\left| \Delta u_m\right|^2\, dx=\lim_{m\to +\infty}\int_{K}|u_m|^{\crito}\, dx =0.
	\end{align*}
\end{step}
\noindent{\it Proof of Step \ref{step:00}:} Note that $u_m \rightharpoonup  0 \hbox{ weakly in } \H$ implies that 
$u_m \to 0\hbox{ strongly in } L^q_{loc}(\e)$ for $1\leq q<\crito$. Therefore, since $2<\crit<\crito$ because $0<s<4$ and the fact $|x|^{-1}$ is bounded on $K$, we get  
\begin{equation}\label{eq:step5.1:0}
	\int_{K}\frac{u_m^{\crit}}{|x|^s}\, dx=o(1) \bb{ and }\int_{K}\frac{u_m^2}{|x|^4}\, dx=o(1)\bb{ as } m\to +\infty.
\end{equation}
We take $\eta \in C^\infty_{c}(\e\backslash\{0\})$ such that $\eta=1$ in $K$ and $0\leq \eta \leq1$. We write $D:=supp(\eta)$. \par 
\medskip\noindent\textit{Step \ref{step:00}.1} We claim $\bb{as } m\to +\infty$ that 
\begin{equation}\label{eq:etaDeltaum}
	\int_{\e}\left| \eta\Delta u_m\right|^2\, dx
	\leq \left( \int_{\e}|\eta u_m|^{\crito}\, dx \right)^{\frac{2}{\crito}}\left( \int_{\e}| u_m|^{\crito}\, dx \right)^{\frac{\crito-2}{\crito}}+\, o(1) .
\end{equation}
\medskip\noindent \textit{Proof of Step \ref{step:00}.1 :} Indeed, using $\lim\limits_{m\to +\infty}\langle \E^{\prime}(u_m), \eta^2u_m\rangle=0$ yields,
\begin{eqnarray}\label{eq:001}
	o(1)&=&\langle \E^{\prime}(u_m), \eta^2u_m\rangle
	=\int_{\e} \langle \Delta u_m, \Delta(\eta^2u_m)\rangle\, dx -\gg\int_{\e}\frac{\left|\eta u_m \right|^2 }{|x|^4}\, dx\nonumber\\
	&&\hspace{4cm}-	\int_{\e}\frac{\eta^2|u_m|^{\crit}}{ |x|^s}\,dx- \int_{\e}\eta^2 |u_m|^{\crito}\, dx.
\end{eqnarray}
Regarding the first term, we have 
\begin{eqnarray}\label{eq:01}
	\int_{\e} \langle \Delta u_m, \Delta(\eta^2u_m)\rangle\, dx&=&\int_{\e}\left| \eta\Delta u_m\right|^2\, dx+\int_{\e} u_m\Delta u_m\Delta(\eta^2)\, dx\nonumber\\
	&&+2\int_{\e}\Delta u_m\langle\nabla (\eta^2), \nabla u_m\rangle\, dx.
\end{eqnarray}
From  H\"older's inequality and since $\nabla u_m \to 0$ in $L_{loc}^2(\e)$ as $m\to +\infty$, we get
\begin{align}\label{eq:02}
	\int_{\e}\Delta u_m\langle\nabla (\eta^2), \nabla u_m\rangle\, dx	&=O\left(\vv \nabla (\eta_m^2)\vv_\infty\int_{supp\left( \nabla(\eta^2)\right)}\left| \Delta u_m\right| \left| \nabla u_m\right| \, dx \right)\nonumber \\
		&=O\left(\vv \nabla (\eta^2)\vv_\infty\vv u_m\vv\left( \int_{supp\left( \nabla(\eta^2)\right) } \left| \nabla u_m\right|^2 \, dx \right)^{\frac{1}{2}}\right)\nonumber\\
		&=o(1) \bb{ as } m\to +\infty.
\end{align}
Also,  since $ u_m \to 0$ in $L_{loc}^2(\e)$ as $m\to +\infty$, we obtain that 
\begin{align}\label{eq:03}
\int_{\e} u_m\Delta u_m\Delta(\eta^2)\, dx=o(1) \bb{ as } m\to +\infty.
\end{align}
Plugging \eqref{eq:02} and \eqref{eq:03} in \eqref{eq:01} yields
\begin{align*}
	\int_{\e} \langle \Delta u_m, \Delta(\eta^2u_m)\rangle\, dx=\int_{\e}\left| \eta\Delta u_m\right|^2\, dx+o(1) \bb{ as } m\to +\infty.
\end{align*}
Hence, it follows from \eqref{eq:001} that 
\begin{equation}\label{eq:prop:1}
	o(1)
	=\int_{\e}\left| \eta\Delta u_m\right|^2\, dx  -\gg\int_{D}\frac{\left|\eta u_m \right|^2 }{|x|^4}\, dx
	-	\int_{D}\frac{\eta^2|u_m|^{\crit}}{ |x|^s}\,dx- \int_{\e}\eta^2 |u_m|^{2^{\star}_0}\, dx.
\end{equation}
Similarly to \eqref{eq:step5.1:0} to get 
\begin{equation*}
	\lim_{m\to +\infty} \int_{D}\frac{\left|\eta u_m \right|^2 }{|x|^4}\, dx=\lim_{m\to +\infty}\int_{D}\frac{\eta^2|u_m|^{\crit}}{ |x|^s}\,dx=0.
\end{equation*}
Therefore, by \eqref{eq:prop:1} and using again the  Hölder's inequality, we find as $m\to +\infty$ that
\begin{align*}
	\int_{\e}\left| \eta\Delta u_m\right|^2\, dx
	&\leq \left( \int_{\e}|\eta u_m|^{\crito}\, dx \right)^{\frac{2}{\crito}}\left( \int_{\e}| u_m|^{\crito}\, dx \right)^{\frac{\crito-2}{\crito}}+o(1).
\end{align*}
 This proves Step \ref{step:00}.1.\qed\par 
\medskip\noindent\textit{Step \ref{step:00}.2} We claim that 
	\begin{equation}\label{eq:Delta(etaum)}
	\int_{\e} \left| \Delta(\eta u_m)\right|^2\, dx =\int_{\e}\left|\eta  \Delta u_m\right|^2\, dx+o(1) \bb{ as } m\to +\infty.
	\end{equation}
\medskip\noindent \textit{Proof of Step \ref{step:00}.2:} Indeed, simple computations yield 
\begin{align}
\int_{\e}\left(  \left| \Delta(\eta u_m)\right|^2-\left|\eta  \Delta u_m\right|^2\right) \, dx  
&=O\left( \int_{\e}\left| \eta \Delta u_m\right| \left|u_m\Delta \eta+2\nabla \eta \nabla u_m \right|\, dx \right.\label{eq:step01:0} \\
&\left. +\int_{\e}\left|u_m\Delta \eta+2\nabla \eta \nabla u_m \right|^2\, dx \right). \nonumber
\end{align}
Using  Hölder's inequality  and  $u_m\to 0$ in $L^2_{loc}(\e)$  and $H^1_{loc}(\e)$ as $m\to +\infty$  
\begin{align}\label{eq:step01:1} 
	\int_{\e}\left| \eta \Delta u_m\right| \left|u_m\Delta \eta+2\nabla \eta \nabla u_m \right|\, dx&=O\left(  \vv \eta \vv_{\infty}\vv \Delta \eta\vv_\infty \vv u_m\vv \left( \int_{D\cap supp\left( \Delta\eta\right)} u_m^2\, dx \right)^{\frac{1}{2}}\right.\nonumber \\
	&\left. +  \vv \eta \vv_{\infty}\vv \nabla \eta\vv_\infty \vv u_m\vv \left( \int_{D\cap supp\left( \nabla\eta\right)} \left| \nabla u_m\right|^2\, dx \right)^{\frac{1}{2}}\right)\nonumber\\
	&=o(1) \bb{ as } m\to +\infty.
\end{align}
In a similar way, we have 
\begin{equation}\label{eq:step01:2}
	\int_{\e}\left|u_m\Delta \eta+2\nabla \eta \nabla u_m \right|^2\, dx=o(1) \bb{ as } m\to +\infty.
\end{equation}
Injecting \eqref{eq:step01:1} and  \eqref{eq:step01:2} in \eqref{eq:step01:0} and we have \eqref{eq:Delta(etaum)}. This proves of  Step \ref{step:00}.2. \qed \par

\medskip\noindent Using \eqref{eq:etaDeltaum} and  \eqref{eq:Delta(etaum)} yield as $m\to +\infty$
\begin{equation}\label{eq:step01:3}
	\int_{\e} \left| \Delta(\eta u_m)\right|^2\, dx
	\leq \left( \int_{\e}|\eta u_m|^{\crito}\, dx \right)^{\frac{2}{\crito}}\left( \int_{\e}| u_m|^{\crito}\, dx \right)^{\frac{\crito-2}{\crito}}+o(1).
\end{equation}
 Now, since $\eta u_m\in \H$, and we go back to the definition of $Q_{\gg,0}(\rr_+^N)$, and $\lim\limits_{m\to +\infty} \int_{D}\frac{\left|\eta u_m \right|^2 }{|x|^4}\, dx=0$, we have that  
\begin{equation}\label{ineq:step:01:sobolev}
	\left( \int_{\e}|\eta u_m|^{\crito}\, dx \right)^{\frac{2}{\crito}}\leq Q_{\gg,0}(\e)^{-1}\vv \eta u_m\vv^2+o(1) \bb{ as } m\to +\infty.
\end{equation}
It follows from \eqref{eq:step01:3} and \eqref{ineq:step:01:sobolev} that
\begin{align}\label{eq:step01:5}
\left[  1-	 Q_{\gg,0}(\e)^{-1} \left( \int_{\e}| u_m|^{\crito}\, dx \right)^{\frac{\crito-2}{\crito}}\right]  \vv \eta u_m\vv^2\leq o(1) \bb{ as } m\to +\infty.
\end{align}
Since $E(u_m)=\beta$ and $E^{\prime}(u_m)=o(1)$ as $m\to +\infty$, we have that 
\begin{align}
	\beta+o(1)&=E(u_m)-\frac{1}{2}\langle E^{\prime}(u_m), u_m\rangle\nonumber\\
	&=\left[\frac{1}{2}-\frac{1}{\crit} \right]\int_{\e}\frac{|u_m|^{\crit}}{|x|^s}\, dx +\left[\frac{1}{2}-\frac{1}{\crito} \right]\int_{\e}|u_m|^{\crito}\, dx.\label{eq:step01:4}
\end{align}
Therefore, since $\crit >2$ when $0\leq s<4$, we obtain as $m\to +\infty$ that 
\begin{align}\label{eq:step01:05}
\int_{\e}\frac{|u_m|^{\crit}}{|x|^s}\, dx\leq 2\, \beta\left[ \frac{N-s}{4-s}\right] +o(1) \bb{ and } \int_{\e}|u_m|^{\crito}\, dx \leq \,  \frac{N}{2} \beta+ o(1).
\end{align}
Therefore, it  follows from \eqref{eq:step01:5} and  $\frac{\crito-2}{\crito}=\frac{4}{N}$  that
\begin{align}\label{eq:step01:6}
	\left[  1-	 Q_{\gg,0}(\e)^{-1} \left( \frac{N}{2} \beta \right)^{\frac{4}{N}}\right]  \vv \eta u_m\vv^2\leq o(1) \bb{ as } m\to +\infty.
\end{align} 
Since $\beta \in (0, \beta^\star)$, and we have $ \left[  1-	 Q_{\gg,0}(\e)^{-1} \left( N\,\frac{\beta}{2} \right)^{\frac{4}{N}}\right]>0$. Moreover, using inequality \eqref{eq:step01:6} yields $\lim\limits_{m\to+\infty} \vv \eta u_m\vv^2=0$. But $\eta=1$ in the compact $K$, then $\lim\limits_{m\to +\infty}\int_k\left| \Delta u_m\right|^2\, dx=0$, from this and by the Sobolev inequality, we obtain that $\lim\limits_{m\to +\infty}\int_{K}|u_m|^{\crito}\, dx =0$. The proof of the Step \eqref{step:00} is complete. \qed\par
\smallskip\noindent For $R>0$, we define 
\begin{align*}
	I_{1,R}:=\lim_{m\to +\infty}\sup \int_{B_R(0)}|u_m|^{\crito}\, dx \bb{ , } I_{2,R}:=\lim_{m\to +\infty}\sup \int_{B_R(0)}\frac{|u_m|^{\crit}}{|x|^s}\, dx,
\end{align*}
and 
\begin{equation*}
	I_{3,R}:=\lim_{m\to +\infty}\sup \int_{B_R(0)}\left(\left|\Delta u_m\right|^2-\gg \frac{|u_m|^2}{|x|^4}  \right)\, dx.
\end{equation*} 
\begin{step}\label{step:01} For $R>0$,  we claim that 
	\begin{align*}
		I_{1,R}^{\frac{2}{\crito}}\leq Q_{\gg,0}(\e)^{-1}I_{3,R} \,\bb{ ;  }\,I_{2,R}^{\frac{2}{\crit}}\leq Q_{\gg,0}(\e)^{-1}I_{3,R}\, \bb{ and }\, I_{3,R}\leq I_{1,R}+I_{2,R}.
	\end{align*}
\end{step} 
 \noindent{\it Proof of Step \ref{step:01}:} Indeed, for $R>0$ we take a cut-off function $\zeta\in C^{\infty}_{c}(\e)$ such that $\zeta=1$ in $B_{R}(0)$ and $0\leq \zeta\leq1$. Since $\zeta u_m \in \H$ and by the definition of $Q_{\gg,0}(\e)$, we get 
 \begin{eqnarray}\label{ineq:step:02:sobolev}
 	\left( \int_{\e}|\zeta u_m|^{\crito}\, dx \right)^{\frac{2}{\crito}}&\leq& 
 	Q_{\gg,0}(\e)^{-1} \left[ \int_{B_R(0)}\left( \left|  \Delta u_m\right| ^2-\gg\frac{| u_m|^2}{|x|^4} \right) \, dx\right. \nonumber\\
 	&&\left. +\, \int_{supp(\zeta)\backslash B_R(0)}\left( \left|\Delta( \zeta u_m)\right| ^2-\gg\frac{| \zeta u_m|^2}{|x|^4}\, dx \right)\right] .
 \end{eqnarray}
It follows from $supp(\zeta)\backslash B_R(0)\subset \e\backslash\{0\}$ and Step \ref{step:00} that
$$\lim_{m\to +\infty} \int_{supp(\zeta)\backslash B_R(0)}\left(  \left|\Delta( \zeta u_m)\right| ^2-\gg\frac{| \zeta u_m|^2}{|x|^4}\right) \, dx =0.$$
Therefore, by \eqref{ineq:step:02:sobolev}  and since $\zeta=1$ in $B_{R}(0)$, we obtain that
\begin{align*}
		\left( \int_{B_R(0)}|u_m|^{\crito}\, dx \right)^{\frac{2}{\crito}}
		&\leq Q_{\gg,0}(\e)^{-1} \int_{B_R(0)}\left( \left| \Delta u_m\right| ^2-\gg\frac{| u_m|^2}{|x|^4} \right) \, dx+o(1),
\end{align*}
as $m\to +\infty$, and we have that $I_{1,R}^{\frac{2}{\crito}}\leq Q_{\gg,0}(\e)^{-1}I_{3,R}$. The proof of $I_{2,R}^{\frac{2}{\crit}}\leq Q_{\gg,0}(\e)^{-1}I_{3,R}$ is similar. \par
\medskip\noindent Since $\zeta^2 u_m\in \H$ and with $\lim\limits_{m\to +\infty}\langle \E^{\prime}(u_m), \zeta^2u_m\rangle=0$, we have 
\begin{align}\label{eq:step:01:01}
	o(1)&=\langle \E^{\prime}(u_m), \zeta^2u_m\rangle\nonumber\\
	&=\int_{\e} \langle \Delta u_m, \Delta(\zeta^2u_m)\rangle\, dx -\gg\int_{\e}\frac{\left|\zeta u_m \right|^2 }{|x|^4}\, dx\\
	&-	\int_{\e}\frac{\zeta^2|u_m|^{\crit}}{ |x|^s}\,dx- \int_{\e}\zeta^2 |u_m|^{\crito}\, dx.\nonumber
\end{align}
It is similar of the proof of Step \ref{step:00}, we have as $m\to +\infty$ that 
$$\int_{\e} \langle \Delta u_m, \Delta(\zeta^2u_m)\rangle\, dx=\int_{\e} \left| \Delta( \zeta^2u_m)\right|^2\, dx+o(1).$$
Therefore, by \eqref{eq:step:01:01}, since $\zeta=1$ in $B_R(0)$, $supp(\zeta)\backslash B_R(0)\subset \rr_+^N\backslash \{0\}$ and with the result of Step \ref{step:00}, we have that 
\begin{align*}
	\int_{B_{R}(0)}\left(  \left| \Delta u_m\right|^2-\gg \frac{\left|\zeta u_m \right|^2 }{|x|^4}\right) \, dx&\leq 
	\int_{\e}\left(  \left| \Delta\left( \zeta^2 u_m\right) \right|^2-\gg \frac{\left|\zeta u_m \right|^2 }{|x|^4}\right) \, dx\\
	&= \int_{B_R(0)}\frac{|u_m|^{\crit}}{ |x|^s}\,dx+ \int_{B_{R}(0)}\left|  u_m\right| ^{\crito}\, dx+o(1).
\end{align*}
Taking $m\to +\infty$ on both sides yields $I_{3,R}\leq I_{1,R}+I_{2,R}$. This proves Step \ref{step:01}.\qed\par
\smallskip\noindent Now, we will to complete the proof of the Proposition \ref{prop2}. Using Step \ref{step:01} yields
\begin{align*}
I_{1,R}^{\frac{2}{\crito}}\leq Q_{\gg,0}(\e)^{-1}I_{1,R}+Q_{\gg,0}(\e)^{-1}I_{2,R},
\end{align*}
this give us
\begin{equation}\label{eq:prop:0}
I_{1,R}^{\frac{2}{\crito}}\left[1-Q_{\gg,0}(\e)^{-1}I_{1,R}^{\frac{4}{N}} \right]\leq Q_{\gg,0}(\e)^{-1}I_{2,R}.
\end{equation}
It follows from \eqref{eq:step01:05} and the definition of $I_{1,R}$ that $I_{1,R}\leq \frac{N}{2}\beta.$ Therefore, by \eqref{eq:prop:0}
\begin{equation*}
	I_{1,R}^{\frac{2}{\crito}}\left[1-Q_{\gg,0}(\e)^{-1}\left( \frac{N}{2}\beta\right) ^{\frac{4}{N}} \right]\leq Q_{\gg,0}(\e)^{-1}I_{2,R}.
\end{equation*}
Since $\beta <\beta^{\star}< \frac{2}{N}Q_{\gg,0}(\e)^{\frac{N}{4}}$, then there exists a constant  $C_1(N,\gg,\beta)>0$ such that $I_{1,R}^{\frac{2}{\crito}}\leq C_1(N,\gg,\beta)\,  I_{2,R}.$ Similar, then there exists $C_2(N,\gg,s,\beta)>0$ such that $I_{2,R}^{\frac{2}{\crit}}\leq C_2(N,\gg,s,\beta)\,  I_{1,R}$. Combining these two inequality we find that 
\begin{equation*}
	I_{2,R} ^{\frac{2}{\crit}}\left[ 1-C_2(N,\gg,s,\beta)C_1(N,\gg,\beta)^{\frac{\crito}{2}}I_{2,R}^{\frac{\crito}{2}-\frac{2}{\crit}}\right]\leq  0.
\end{equation*}
Therefore, since $\crito>\frac{4}{\crit}$ we have $I_{2,R}=0 \bb{ or there exists } \ep:=\ep(N,\gg,s,\beta) \bb{ such that } I_{2,R}\geq \ep.$
Similarly, we have $I_{1,R}=0 \bb{ or there exists } \ep \bb{ such that } I_{1,R}\geq \ep.$ This ends of the proof of Propostion \ref{prop2}.\qed\par 
\medskip \noindent \underline{\textbf{End of proof of Theorem \ref{theo2}} :} Indeed, we let $(u_m)_{m\in \nn}$ be the Palais-Smale sequence for $\E $ that was constructed in Proposition \ref{prop0}. First, we claim that 
\begin{eqnarray}\label{eq:crit000}
	\lim_{m\to +\infty} \sup\int_{\e}\left|u_m \right|^{\crito}\, dx >0. 
\end{eqnarray}
	 \smallskip\noindent Indeed, otherwise $\lim\limits_{m\to +\infty} \sup\int_{\e}\left|u_m \right|^{\crito}\, dx =0.$ Using again $\lim\limits_{m\to +\infty}\langle \E^{\prime}(u_m), u_m\rangle=0$ yields,
	\begin{equation*}
		\int_{\e}\left(  \left| \Delta u_m\right|^2-\gg \frac{\left|u_m \right|^2 }{|x|^4}\right) \, dx= \int_{\e}\frac{|u_m|^{\crit}}{ |x|^s}\,dx+o(1).
	\end{equation*}
Therefore, we go back to the definition of $Q_{\gg,s}(\e)$
\begin{align*}
\left( \int_{\e}\frac{|u_m|^{\crit}}{ |x|^s}\,dx\right)^{\frac{2}{\crit}}  
	\leq  Q_{\gg,s}(\e)^{-1} \int_{\e}\frac{|u_m|^{\crit}}{ |x|^s}\,dx+o(1).
\end{align*} 
This give us, 
\begin{equation*}
\left( \int_{\e}\frac{|u_m|^{\crit}}{ |x|^s}\,dx\right)^{\frac{2}{\crit}}\left[1-Q_{\gg,s}(\e)^{-1}	\left( \int_{\e}\frac{|u_m|^{\crit}}{ |x|^s}\,dx\right)^{\frac{\crit-2}{\crit}} \right] \leq o(1).
\end{equation*}
It follows then from the left inequality of \eqref{eq:step01:05} that
\begin{equation*}
\left( \int_{\e}\frac{|u_m|^{\crit}}{ |x|^s}\,dx\right)^{\frac{2}{\crit}}	\left[1-Q_{\gg,s}(\e)^{-1}\left(2 \frac{N-s}{4-s}\beta\right) ^{\frac{\crit-2}{\crit}} \right]\leq o(1). 
\end{equation*}
Since $0<\beta<\beta^{\star}$, we have the quantity between the brackets is positive. Thus, we get $	\lim\limits_{m\to +\infty} \int_{\e}\frac{|u_m|^{\crit}}{ |x|^s}\,dx=0.$ Therefore, using \eqref{eq:step01:4} yields  $\beta=0$ which contradicts the fact that $\beta \in (0,\beta^{\star})$. This proves the claim.\qed\par
\medskip\noindent Next, we claim that the sequence $(u_m)_{m\in \nn}$ is bounded in $\H$. \par
\smallskip \noindent Indeed, since $u_m$ is a Palais-Smale sequence for $\E$ and using \eqref{eq:coercive} yields
\begin{align*}
	\beta+o(1)&=\E(u_m)-\frac{1}{\crit}\langle \E^{\prime}(u_m), u_m\rangle\nonumber\\
	&=\left[\frac{1}{2}-\frac{1}{\crit} \right]\int_{\e}\left( \left| \Delta u_m\right|^2-\gg\frac{u_m^2}{|x|^4}\right)\, dx +\left[\frac{1}{\crit}-\frac{1}{\crito} \right]\int_{\e}|u_m|^{\crito}\, dx\\
	&\geq C \left[\frac{1}{2}-\frac{1}{\crit} \right]\int_{\e}\left| \Delta u_m\right|^2\, dx,
\end{align*}
where $C$ is a positive constant. It follows then from $2<\crit<\crito$ that $u_m$ is bounded in $\H$. This proves the claim. \qed\par 
\medskip \noindent Since $u_m$ is bounded in $\H$, then there exists $u\in \H$ such that $u_m \rightharpoonup  u\hbox{ weakly in } \H$. If $u\not\equiv 0$, we get that $u$ is  a nontrivial weak solution of \eqref{eq:doublehardybi}. \par 
\smallskip\noindent If $u\equiv0$, we have $u_m \rightharpoonup  0\hbox{ weakly in } \H$. We claim that, for small enough $\ep^\prime>0$, there exists another Palais-Smale sequence $(v_m)_{m\in \nn}$ satisfying the properties of Proposition \ref{prop2} and 
$$\int_{B_1(0)}\left|v_m \right|^{\crito}\, dx =\ep^{\prime} \bb{ ; $v_m$ is bounded in $\H$ for all } m\in \nn.$$ 
\noindent Indeed, by \eqref{eq:crit000}, we can take $c:=\lim\limits_{m\to +\infty} \sup\int_{\e}\left|u_m \right|^{\crito}\, dx$. We set $\ep_0:=\min\{c, \frac{\ep}{2}\}$, where $\ep>0$ is the same which we obtain from Proposition \ref{prop2}. Therefore, for any $\ep^\prime \in (0, \ep_0)$, there exists a sequence $(r_m)_{m\in \nn}>0$ such that up to a subsequence $\int_{B_{r_m}(0)}\left|u_m \right|^{\crito}\, dx =\ep^\prime$. Define now $v_m(x):=r_m^{\frac{N-4}{2}}u_m(r_mx) \bb{ for all } x\in \rr_+^N.$ With change of variable, we write 
\begin{equation}\label{eq:0}
\int_{B_{1}(0)}\left|v_m \right|^{\crito}\, dx =	\int_{B_{r_m}(0)}\left|u_m \right|^{\crito}\, dx =\ep^\prime.
\end{equation}
As one checks, $(v_m)_{m\in \nn}$ is also a Palais-Smale sequence for $\E$ that satisfies the properties of Proposition \ref{step:00}. Using the definition of $v_m$ and the boundedness of the sequence $u_m$ yields $v_m$ is bounded in $\H$. This ends the prove of Claim. \qed \par 
 \smallskip \noindent Hence, we can assume that there exists $v\in \H$ such that, up to a subsquence $v_m \rightharpoonup  v\hbox{ weakly in } \H$.\par
\smallskip \noindent We claim now that $v$ is a nontrivial weak solution of problem \eqref{eq:doublehardybi}.\par
\smallskip\noindent Indeed, if $v\equiv 0$. It follows from the result of Proposition \ref{prop2} that 
\begin{eqnarray*}
	\hbox{ either }\lim_{m\to +\infty}\sup \int_{B_1(0)} \left| v_m\right|^{\crito}\, dx  =0 	\hbox{ or }	\lim_{m\to +\infty}\sup \int_{B_1(0)} \left| v_m\right|^{\crito}\, dx \geq \ep.
\end{eqnarray*}
Since $\ep^\prime \in (0, \frac{\ep}{2})$, this is contradiction with \eqref{eq:0}. Then $v\not\equiv 0$.\qed\par 
\smallskip \noindent Since $(v_m)_{m\in \nn}$ is a sequence Palais-Smale for $\E$, we have 
\begin{eqnarray}
o(1)&=& \langle \E^\prime(v_m), \varphi \rangle \nonumber\\
&=&\int_{\e} \langle \Delta v_m, \Delta \varphi \rangle \, dx-\gg\int_{\e}\frac{v_m\varphi}{|x|^4}\, dx\label{eq:000}\\
&&-\int_{\e} \left|v_m \right|^{\crito-2}v_m\varphi\, dx-\int_{\e}\frac{\left|v_m\right|^{\crit-2}v_m\varphi}{ |x|^s}\,dx,\nonumber
\end{eqnarray}
for all $\varphi\in C^\infty_c(\e)$. Using $v_m \rightharpoonup  v\hbox{ weakly in } \H$ yields
\begin{equation}\label{eq:convergefaible}
\lim_{m\to +\infty}	\int_{\e} \langle \Delta v_m, \Delta \varphi \rangle \, dx = \int_{\e} \langle \Delta v, \Delta \varphi \rangle \, dx \bb{ for all }\varphi\in C^\infty_c(\e).
\end{equation}
Since $v_m$ is bounded in $\H$, we get that $v_m$, $\left|v_m \right|^{\crito-2}v_m$ and $\left|v_m\right|^{\crit-2}v_m$ are bounded in  $L^2( \e,|x|^{-4})$, $L^{\frac{\crito}{\crito-1}}(\e)$ and $L^{\frac{\crit}{\crit-1}}( \e,|x|^{-s})$ respectively. Therefore, we get that 
\begin{eqnarray*}
	\left\{\begin{array}{ll}
		v_m \rightharpoonup u \bb{ weakly in } L^2( \e,|x|^{-4}),\\
		\left|v_m \right|^{\crito-2}v_m \to\left|v \right|^{\crito-2}v \bb{ weakly in } L^{\frac{\crito}{\crito-1}}(\e),\\
		\left|v_m\right|^{\crit-2}v_m\rightharpoonup\left|v\right|^{\crit-2}v \bb{ weakly in }
		L^{\frac{\crit}{\crit-1}}(\e,|x|^{-s}).
	\end{array}\right.
\end{eqnarray*}
Moreover, passing the $m\to +\infty$ in \eqref{eq:000} and using \eqref{eq:convergefaible} yields
\begin{eqnarray*}
	o(1)
	&=&\int_{\e} \langle \Delta v, \Delta \varphi \rangle \, dx-\gg\int_{\e}\frac{v\varphi}{|x|^4}\, dx\\
	&&-\int_{\e} \left|v \right|^{\crito-2}v\varphi\, dx-\int_{\e}\frac{\left|v\right|^{\crit-2}v\varphi}{ |x|^s}\,dx \bb{ for all }\varphi\in C^\infty_c(\e).
\end{eqnarray*}
Thus, $v$ is a weak solution of \eqref{eq:doublehardybi}. This completes the proof of Theorem \ref{theo2}.\qed
\section{Appendix}
\begin{prop}\label{prop:nonext}
	Let $\Omega \subset\rr^N$ be a smooth domain such that $0\in \partial \Omega\neq\emptyset$ (No bound-
	edness is assumed). If $\gg \leq 0$ and $s=0$, we have that $Q_{\gg,0}(\Omega)= S_N$, where $S_N$ is defined in \eqref{eq:bestconstantsoblevrn} and there is no extremal.
\end{prop}
\medskip\noindent \textit{Proof of Proposition \ref{prop:nonext}:} Indeed, we call back  $\crito=\frac{2N}{N-4}$. Since $\gg\leq 0$, we
find 
\begin{eqnarray*}
	\frac{\int_{\Omega} \left( |\Delta u|^2-\gg \frac{u^2}{|x|^4}\right) \, dx}{\left(\int_\Omega |u|^{\crito}\, dx \right)^{\frac{2}{\crito}}}\geq \frac{\int_{\Omega}  |\Delta u|^2 \, dx}{\left(\int_\Omega |u|^{\crito}\, dx \right)^{\frac{2}{\crito}}}\geq S_N,
\end{eqnarray*} 
and we have that $Q_{\gg,0}(\Omega)\geq S_N$. Fix $x_0\in \Omega$ such that $x_0\neq 0$. We define $(U_\ep)_\ep$ as in \eqref{def:Ue}. It follows from \eqref{eq:i1ep} that 
$\lim\limits_{\ep\to 0}\int_\Omega \frac{U_\ep^2}{|x|^4}\, dx=0$. Moreover, with \eqref{eq:la} and \eqref{eq:U2**}, we get that
\begin{equation*}
\lim_{\ep\to 0} \frac{\int_{\Omega}|\Delta U_\ep|^2\, dx}{\left(\int_\Omega |U_\ep|^{\crito}\, dx \right)^{\frac{2}{\crito}} }=S_N.
\end{equation*} 
We then get $Q_{\gg,0}(\Omega)\leq S_N$. This proves that $Q_{\gg,0}(\Omega)=S_N$. If there was an extremal for $Q_{\gg,0}(\Omega)$, it would
also be a extremal for $S_N$, with their support is the whole of $\rr^N$, contradicting since the extremal has support in $\Omega\neq \rr^N$ smooth. This proves Proposition \ref{prop:nonext}.\qed

\bibliographystyle{abbrv}

\end{document}